# Floating-Point Numbers
# with Error Estimates
# (revised)

Glauco MASOTTI

**Abstract** — The original work [25] is here reconsidered, so many years after. The old text has been revised, plus several considerations have been added, in order to clarify some controversial aspects of the work, and to envision possible developments. A section has been added, to review the effects of the original paper.

The study addresses the problem of precision in floating-point (FP) computations. A method for estimating the errors which affect intermediate and final results is proposed and a summary of many software simulations is discussed. The basic idea consists of representing FP numbers by means of a data structure collecting value and estimated error information. Under certain constraints, the estimate of the absolute error is accurate and has a compact statistical distribution. By monitoring the estimated relative error during a computation (an ad-hoc definition of relative error has been used), the validity of results can be ensured. The error estimate enables the implementation of robust algorithms, and the detection of ill-conditioned problems. A dynamic extension of number precision, under the control of error estimates, is advocated, in order to compute results within given error bounds. A reduced time penalty could be achieved by a specialized FP processor. The realization of a hardwired processor incorporating the method, with current technology, should not be anymore a problem and would make the practical adoption of the method feasible for most applications.

**Preamble**

What has brought me to revise my old paper [25], more than 18 years after, and publish this document? A part some casual circumstances, the main reasons are perhaps the conviction that this is still a good idea, and also my sufferance in seeing my past work so much neglected or misunderstood in all these years.

This document is thus also a prayer to experts in the field and to hardware architects to reconsider this work for practical implementation.

My original paper had a hard time since the beginning and it was refused for publication a couple of times, before encountering the attentive consideration of Prof. Tony C. Woo, who cured the special issue of CAD on *"Uncertainties in geometric design",* and decided to publish it, in spite of a strong criticism from some of the peer reviewers, that I had to counter[1]. I am going to report in footnotes some of their remarks and my replies to them. In my humble opinion, it appears evident that most of the criticism was inconsistent, and that my work was poorly understood. This shows once more all the limits of the peer review process [28]. Months or years of work risk often to be slashed by some reviewers who read your paper without attention, perhaps in their spare time, and don't understand the essence of your work. The risk is higher if the ideas that you expose and the language and formalism you use are new and non-orthodox, conditions which may apply to this work. There are no filters here, so anyone of you can judge if this is true or not and the value of this work.

The original work [25] is here reconsidered. The old text has been revised, plus several considerations have been added, in order to clarify some controversial aspects of the work — as they appeared, prior to publication, in the criticism of the peer review and later in some citations — and also to envision possible developments. A section has been added, to review the effects of the original paper. Finally the concluding remarks have been updated, taking into account all this.

From now on, like in the original paper, I will switch to using the formal and more elegant "we", in place of the more realistic, but not nice, "I". Nevertheless I found sometimes more appropriate the use of "I", like in this preamble and in most of the footnotes.

## 1. State of the art[2] and introduction

The errors that unavoidably affect floating-point (FP) computations are a well known source of troubles for all numerical algorithms [1], [2], [3]. Studies on systematic and statistical properties of various FP systems [4], [5], led to the IEEE FP standard [6], where the best possible strategies have been adopted to reduce the amount of rounding errors.

Coping with roundoff errors is a general problem, but some applications suffer more than others from the consequences of imprecise numerical computations. This is the case, for instance, of geometric applications [7] (from where most of the author's experience derives), because in this field we are processing problems in a continuous space with finite precision, while this kind of problems would require infinite precision to be solved properly. Anyone who has been involved in the implementation of geometric algorithms, in the fields of computer graphics, CAD-CAM, simulation and modeling, knows the frustration induced by unexpected program faults due to numerical errors, and the timeless effort required in trimming an application to work in all possible conditions.

---

[1] This caused a quite delayed publication of the results of my work, which in reality dates back to 1991 [38], however the criticism has at least stimulated improvements in the presentation of the work.
[2] The state of the art is not changed much since 1993. The main techniques used in practice are substantially the same, so I made here little variations with respect to what I wrote in the original paper. For those interested, a somewhat more comprehensive and updated survey of known methods, but mainly focused on geometric applications, can be found in [33].

Roundoff errors, originated at a certain point in a computation, propagate to subsequent steps. They are amplified in ill-conditioned problems [1], [7] up to the point of making these problems intractable.

We usually know that numerical errors are affecting the results, but we do not know how big these errors are. The idea of resorting to *error estimates* is not new. The concept of robustness of numerical algorithms is connected to *verification* of correctness of the computed output. It is diffusely perceived that if we could know, at important steps of a computation, how many significant digits are in the computed values, or, more precisely, the amount of the error accumulated in the results, we could be able to write robust algorithms. However, practical methods, which could be used easily and extensively in applications, to obtain accurate estimates of these errors, have not been hitherto available.

One of the first methods for error estimation was the so called *backward analysis* proposed by Wilkinson [8]. Moore [9] later introduced the concept of *interval arithmetic*. These methods provide upper bounds for the accumulated roundoff error. Their results are thus too pessimistic, as they do not account for error compensations. The output of a computation usually satisfies much tighter error bounds [5]. Moreover the implementations of these methods [10], [11] are quite time consuming. Robustness cannot in any case be achieved for free and, as reported in the literature [12], there are cases where the implementation of robust algorithms costs as much as a factor of 100 in speed.

Vignes and La Porte [13] introduced the *permutation-perturbation* method which appears as a major advancement in error analysis. The method can estimate the number of significant digits in the result of any algorithm [14]. The quantification of the total cost of applying the method is not explicitly provided in the referenced papers, but estimating the number of significant digits with a confidence level of 0.95 appears to introduce a time penalty factor *TPF* of at least 3 in total execution time. The method in fact requires re-executing three times the computations of an algorithm, randomly "perturbing" the least significant digit of the mantissa of FP numbers and the order of execution of operations. Taking into account the overhead of control, it follows that we should have *TPF* > 3, for a 0.95 confidence level. Moreover a 0.95 confidence level for the error estimation is, most likely, not sufficient for many applications. "Absolute" robustness requires confidence levels very close to 1, to reduce the probability of failure to a negligible level. We deem that, in order to achieve such very high confidence levels, the number of samples which have to be evaluated would have to be increased considerably, proportionally increasing the cost of the method.

The previous considerations state that, in the current situation, it is too expensive to compute error estimates by inserting appropriate steps throughout the code of a program, either in terms of cost for setting up the software, or of time penalties which are introduced, or both. It follows that computing error estimates is impractical and thus usually not done in applications.

Nevertheless we need error estimates at certain points. We are thus forced to guess gross upper bound for numerical errors, but it is almost a blind guess.

The correctness of an algorithm is mainly endangered when the time comes to make critical *decisions* on the basis of imprecise data. As we made somewhat arbitrary assumptions on error estimates, we are running the risk of making the wrong choice. Let us take an example. The typical case is when we must *compare* two FP numbers, which result from computations of unknown condition numbers [1], [7], and thus affected by unknown errors. Problems which present this situation are for instance: comparing two points in space for coincidence (using a distance between them); testing for inclusion of a point in a set (e.g. deciding whether a point lies over a surface); testing for parallelism or orthogonality of curves and surfaces, etc. This kind of decisions may be critical for the correctness of an algorithm as completely different results could be produced making one choice or another. The point is that *discrete* changes in topology or other geometric characteristics are induced by continuous changes of data.

What is customarily done in these circumstances is to introduce some sort of tolerance, usually making the supposition that a combination of the relative and absolute errors were less than a fixed small number $\varepsilon$, which we are used to call *"epsilon"*. Combining absolute and relative error is common programming practice (a glance at the code of some geometric modeling systems, e.g. PADL-2 [15], confirms this statement), as relative errors do not work in a neighborhood of the value 0. Let us indicate with $x$ and $y$, two FP numbers, and with **x**,**y**, the true real values that they represent. The previous argument is to say that we can consider:

$$\mathbf{x} = \mathbf{y} \quad \Leftrightarrow \quad |x - y| < \varepsilon(1 + |x|) \tag{1.1}$$

for some number $\varepsilon$ which is sufficiently *small*, but *not too small*. Here lies the pitfall: how big should $\varepsilon$ be? If we make it too big it is like considering the two numbers less accurate than they actually are, and we are running the risk of considering them as equal, while they identify two different entities. On the other hand, if we make $\varepsilon$ too small, sometimes the error may exceed the tolerance in (1.1), making us fail to consider the two numbers as equal when they should identify the same entity. In practice we choose $\varepsilon$ by experience, as a compromise, trying to cope with all the situations that may arise in practice. The consequences are that even after a fine-tuning and a lot of time spent with several test trials, we cannot guarantee the correct behavior of an algorithm in all operating conditions, i.e. in certain critical situations our program may *occasionally fail*. This is particularly true if the problem at hand is numerically ill-conditioned, and, even worse, we cannot even detect ill-conditioned problems in general. This situation is described efficaciously by Dobkin and Silver [12], as follows:

> "The most widely applied solution to the problem of roundoff error is the ad hoc approach: calling the local guru[3] to pull a fix out of his/her magic box. This usually entails arbitrary increasing precision, re-ordering calculations, tweaking specific numbers, or arbitrarily selecting *epsilon* values, and in most instances, will only solve a set of problems temporarily and does not attack the underlying cause of the roundoff error. Needless to say, this approach is far from robust and consistent."

This is to say that the problem has been tackled hitherto mostly by means of heuristic and empirical methods, based on programming practice and, more or less, sound judgment. It has been left almost untouched by theoretical developments and no new technology, suitable for general use, and capable of eliminating or greatly reducing the effects of imprecise computations, has been developed.

The most simplistic approach to cope with inaccuracy of computations is to increase the precision of representation using *larger formats* for FP numbers [36], [37]. Today computers, however, do not efficiently support precision extension beyond a certain format, so that critical computations are customarily done at maximum available precision. Proposing the *indiscriminate* use of larger and larger FP formats, however, besides its practical cost, does not face the problem, it attempts to avoid it, but the problem would reappear a little further forward. Increasing the size of the mantissa, in fact, does not provide us error estimation, so we cannot know what the available precision of computed values is, no matter how large a format we use. The practical result is that we are simply narrowing down the traps which will make our program fail, but we will still be unaware of the amount of the errors and we will have no control on the phenomenon of error generation and propagation.

Accuracy and robustness issues, particularly in geometric computations, have received attention by researchers. A number of approaches have been proposed, among which the use of interval arithmetic [16], rational numbers [17] and purely symbolic computations [7], [29], which are alternatives to the employment of FP numbers. An in depth review of all these studies is beyond the scope of this paper. Known approaches are reviewed in [12], and more extensively presented and discussed by Hoffmann in [7] and [18]. We can briefly summarize the results as follows:

— Interval arithmetic tends to produce overly pessimistic results, because, as we said, intervals account for the worst possible case.

---

[3] In fact I used to be such a "local guru" for years, so I know this situation very well, I fully subscribe the citation.

- Rationals cannot (by definition) represent precisely transcendental numbers. However, even restricting the domain of initial representable data to rationals with a bounded denominator, they require an astoundingly high number of bits to perform a computation exactly, otherwise they do not avoid roundoff problems. This is in practice equivalent to performing FP computations with unbounded precision carrying along all the computed bits throughout any computation. This is thus too expensive in terms of memory usage and computation time.
- Symbolic reasoning helps in increasing robustness, but it can be used only to a limited extent. Avoiding completely numerical computations seems unfeasible, as this raises complex theoretical and computational problems [18].
- Other solutions like data normalization [19], computing with uncertainty regions [20], etc. [33], look like formalizations of ad-hoc approaches. They reported partial success in specific domains, without guaranteeing absolute robustness, and are not for general use in numeric computations.

With the method proposed here instead we can get *accurate estimates* of the absolute error which affects the result of a FP computation. A tight interval for the error can be evaluated with a confidence level which can closely approach 1. This can be achieved with minimal changes to existing software or to current methodology in writing new software. Knowing the amount of the error, we can avoid such approximations as in (1.1). We can write robust algorithms, being certain of the validity of the decisions that have to be taken, or otherwise *detect* ill-conditioned problems, which require more precision in the representation of numbers to be solved properly.

This is the case when the uncertainty on the computed values invalidates their meaningfulness. The estimated error, in these situations, can thus serve to drive *dynamic extension* of precision, as required by the problem at hand, in order to compute final results within given error bounds. This strategy is proposed and analyzed also in [12]. The error estimation in [12] is based on Vignes's method [13], but the observations concerning dynamic precision extension maintain their validity even if a different error estimation method, like ours, is used. Researches on hardware implementations for variable precision arithmetic have been reported already many years ago [21], [22], as well as in recent times [35].

Our method has been tested with a software simulation over a number of cases, getting encouraging results. To make it usable in practice however, a fast implementation should be feasible. This is our case. In fact, it turns out that the method is suitable for hardware implementations of various levels of complexity and performance, which could substantially reduce the penalty factor, with respect to performing standard FP operations, down to a negligible or acceptable level for most applications.

## 2. Error estimates

Let us suppose that, at a certain point in our program, as the result of a certain computation, we have a variable with a FP value *x,* while the true exact value of the computation should be a real number **x**. If we know that this value is affected by an absolute error **e**, we have: **x** = *x* + **e**. In practice however only an approximation of **e** is possible, because again we have roundoff errors in it. Let us call *e* this FP representation of the error. The couple:

$$(x, e) \tag{2.1}$$

contains this information. Let us suppose from now on that we have $T$ bits available for representing the mantissa of $x$ and $T_e$ bits available for the mantissa of $e$.

When the FP value *x* is used in a FP computation, the error which affects *x* propagates to the result of the computation. The propagation of roundoff errors can be evaluated by applying a simple differential analysis [1]. In particular if two FP numbers *x* and *y,* affected respectively by the errors $e_x$ and $e_y$, are operands of FP arithmetic operations, their errors propagate to the result with the contribution *pe,* given by the following formulas:

$$pe(x+y) \;=\; e_x + e_y \tag{2.2}$$

$$pe(x-y) \;=\; e_x - e_y \tag{2.3}$$

$$pe(x*y) \;\approx\; y*e_x + x*e_y \tag{2.4}$$

$$pe(x/y) \;\approx\; (e_x - (x/y)*e_y)/y \tag{2.5}$$

and for the square root of *x* we have:

$$pe(\sqrt{x}) \;\approx\; (1/(2*\sqrt{x}))*e_x \tag{2.6}$$

and so on for other functions. Note that equations (2.2)…(2.6) are for the *absolute* errors, while textbooks on numerical analysis usually report formulas for the *relative* error. The reasons for using absolute errors in our reasoning will appear clear later on. The above formulas are exact for addition and subtraction, and first order approximations for the other operations. To the propagated error *pe* we have to add the local roundoff error *le,* produced to store the result of this step. If **z** is the exact result of a real operation, the local roundoff error *le* of the correspondent FP operation is given by:

$$le = \mathbf{z} - rd(\mathbf{z}) \tag{2.7}$$

where *rd()* is a FP rounding function defined over the domain of the real numbers **R**. Therefore the total estimated error in the result *z* of each operation is given by:

$$e_z = pe + le \qquad (2.8)$$

Equations (2.3)…(2.8) tell us how to update the error estimate throughout every step of a computation.

## 3. Determination of local roundoff error

Equation (2.7) defines *le* in terms of the exact result *z* of a computation, which is unavailable. How can we skip this obstacle? Here we have to introduce one of the basic hypotheses of our method. We are assuming that, in general, we can *temporarily* have at our disposal the result *w* of a FP operation with a number of bits $T_w$ significantly greater than *T*. It could not be possible in certain cases to satisfy this condition, but it should be verified in particular when a significant roundoff error is about to be made. We know that a multiplication of two FP numbers *x* and *y*, with a mantissa of *T* bits, typically produces a result with 2\**T* bits. A division can produce an infinite sequence of bits. Addition and subtraction produce a result with $T+|exp_x - exp_y|$ bits, where $exp_x$, $exp_y$ are the exponents of *x* and *y*. What is a reasonable choice for $T_w$? As a reasonable compromise, we can assume that the result of an arithmetic operation is temporarily available in an accumulator with a mantissa of 2\**T* bits. The components of the result of a computation, according to format (2.1) are then given by:

$$z = rd(w) \qquad (3.1)$$
$$le = z - w \qquad (3.2)$$
$$(z, rd_e(pe + le)) \qquad (3.3)$$

where *rd()* is a rounding function which operates on the accumulator *w* and produces a result *z* with a mantissa of *T* bits, $rd_e()$ is a rounding function which produces a result with a mantissa of $T_e$ bits. In the hypothesis $T_w = 2*T$, *le* can be represented with up to *T* bits, so we can assume $T_e \leq T$.

It should be pointed out that, in order to determine the estimated roundoff local error *le*, we do not really need to compute the difference *z – w*. In fact, let us suppose that the mantissa $m_w$ of the accumulator *w* is composed of two parts *u* and *v*, each one *T*-bit long:

$$m_w = [u|v] \qquad (3.4)$$

where: $u = 0.u_1 u_2 \ldots u_T$ and $v = 0.v_1 v_2 \ldots v_T$. Recalling the definition of a rounding function [1], we have:

$$m_z = \begin{cases} u & \text{if } v_1 = 0 \\ u + 2^{-T} & \text{if } v_1 = 1 \end{cases}$$

and consequently, the (signed) mantissa of *le* is:

$$m_{le} = \begin{cases} v & \text{if } v_1 = 0 \\ v - 2^{-T} & \text{if } v_1 = 1 \end{cases}$$

which are very fast operations (note that $m_{le}$ changes sign in the second case).

We have not considered yet the problem of initializing the error estimate at the beginning of each computation sequence. At the beginning we necessarily get values from constants or from input data. Therefore the internal representation of numerical data should use the same format (2.1) used to represent values of variables and intermediate results. In order to obtain an estimate of the initial conversion error, the compiler (in the case of constants) or the input handler (in the case of input data) should provide a preliminary representation of the data with a number of bits exceeding those representable in the mantissa of the FP format used. The best conditions arise if we have available a representation for the initial values, with at least $T+T_e+1$ bits, so we can initialize both *x* and *e* with the necessary accuracy. Let us call the value of this extended temporary representation *w*. Under the previous assumptions we have enough information to initialize the fields in (2.1) according to (3.1)…(3.3), setting *pe* = 0.

**4. Accuracy of the error estimation**

The error estimation gives us only an approximation of the true error associated with a FP number. In fact, it is, in its turn, affected by errors due to the first order approximation used in error propagation formulas and to roundoff errors in its evaluation.

The effect of second and higher order terms and of roundoff errors over the estimated error can be taken into account *only statistically* (as we will see later). It is in this way that the behavior of error estimation has been studied in many test cases.

After an ill-conditioned calculation or a very long sequence of computations the error estimation can accumulate these errors to an extent that makes it meaningless, just in the same way as it happens for the value of the FP number with which the error estimation is associated. How can we keep this phenomenon under control?

From now on, for clarity, let us call the estimated error *ee,* as opposed to the true error *e.* First of all, for the validity of (2.2)…(2.6), in order to have *e* ≈ *ee*, the estimated errors should be *small*. What does "small" mean? The meaning of small depends on the value of *x.* It seems reasonable to require that the estimated relative error of *x:*

$$re(x, ee) = |ee/x|, x \neq 0 \qquad (4.1)$$

be small, but with the definition (4.1) we are in trouble if *x* = 0. We should also expect a fuzzy behavior as *x* → 0. In fact, if *x* ≈ 0, the question is whether this is a corrupted representation of the real value **0** or whether it represents a small real number. Let us examine more closely what happens in the following cases:

1) |*x*| » 0. In this case we can assume *re* = |*ee/x*| without any problem, and require that *re* < *RTHD* « 1. As long as this condition holds we should also have *x* ≈ **x**. The constant *RTHD* represents a threshold of acceptability for *re.* It should be small enough to be confident about the results we have, but not too small in relation with the available precision. Its value is quite independent from the application; it depends mainly on the precision used to represent the FP numbers, i.e. *T* and $T_e$ in particular (later on we will discuss the implications of varying *RTHD*).

2) *x* = 0, **x** = 0. In this case the definition (4.1) cannot be applied. Here we necessarily have to use the estimated absolute error *ee* as a measure of accuracy of *x.* Sometimes we may have concomitantly *x* = 0, *ee* = 0, but we may have *x* = 0, *ee* ≠ 0 as well. In order to be confident about the calculated value *x* = 0, we require that also the estimated error be very small, i.e. *ee* < *EPS,* where *EPS* is so small that we can consider a smaller value as equivalent to 0. Therefore *EPS* depends to some extent on the applications and the range of meaningful values of *x.* If it can be set many orders of magnitude smaller than any meaningful value of *x,* this dependence is negligible. This is the case when the precision used in our FP representation is high enough.

3) *x* ≈ 0, **x** = 0. In this case the true value **x** of *x* is 0, but *x* ≠ 0. How can we detect this situation? The decisive point here is that as long as *ee* is correct, *x* and *ee* tend to assume *symmetric* values with respect to 0, i.e. if *x* → 0, *ee* → -*x*, hence the relative error |*ee/x*| → 1, and so we cannot use it. However, for this same reason, we have |*x+ee*| « |*ee*| ≈ 0, so it seems reasonable in these

conditions to require $|x+ee| < EPS$, as a measure of accuracy. If we want to achieve uniformity in the definition of *re*, so that we can require *re* < *RTHD* also in this case, we should define *re* = $|x+ee|/EEZ$, where *EEZ* = *EPS*/*RTHD*. It should be pointed out here that having exactly $x = -ee$, which would imply *re* = 0, is a very unlikely event — in fact it never occurred in our simulations (see the observed values of *re* in Figs. 3÷8) — as the two numbers are computed throughout two independent sequences of computations. This helps. In fact, if in a long sequence of computations the value of *x* increases, due to the accumulation of errors, the difference between the moduli of *x* and *ee* is amplified by the factor 1/*EEZ*, making *re* exceed the threshold *RTHD* at a certain point, i.e. when $||x|-|ee|| > EPS$.

4) $x \approx 0$, **x** ≠ 0.   Here the true value of *x* is a small real number **x** ≠ 0, hence *x* and *ee* should not have a symmetric disposition around 0. We should have $|ee| \ll |x|$, even if $x \approx 0$, so we can apply the common definition (4.1) of *re*.

5) $x = 0$, **x** ≠ 0.   The true value of *x* is a small real number **x** ≠ 0, as in the previous case, but here the computed value happens to be 0. This is a real coincidence, so that this case should be very unlikely, in fact it never occurred in our simulations. In any case we would have here $|ee/x| = \infty$ (unless *ee* = 0 too, which means the concurrence of two improbable coincidences!), but a measure of the relative error can be provided by the alternative definition used at point 3). In the case at hand we have *re* = $|x+ee|/EEZ$ = $|ee|/EEZ$. Requiring *re* < *RTHD*, in this case too, implies $|ee| < EPS$, so that, as long as $e \approx ee$, it is $|x| < EPS$. This means that it should not be possible to have this case for **x** ≫ 0.

Summing up we should use two different formulas for computing *re* :

$$re(x, ee) = |ee/x|; \quad x \neq 0, \mathbf{x} \neq 0$$

or:

$$re(x,ee) = |x+ee|/EEZ; \quad x = 0 \mid (x \approx 0, \mathbf{x} = 0)$$

according to the case at hand, which also depends on the true value **x**, which we do not know! Nevertheless we can resolve the ambiguity by choosing the smallest value given by the two formulas. In fact it can be verified, on the basis of the discussion above, that when one of the two formulas should be used, the other one gives a larger value. We can thus assume the following generalized definition of *re(x,ee)* as being valid in any case:

$$re(x,ee) = min(|ee/x|, |x+ee|/EEZ) \qquad (4.2)$$

where $EEZ = EPS/RTHD$. Equation (4.2)[4] implies that we do not need to distinguish the various cases in practice, in order to compute *re*! This fact is of fundamental importance because we would be unable to make the distinction, as the true value **x** is unknown.

Let us suppose at this point that in the course of each step *i* of a computation the following relation holds:

$$re_i = re(x_i, ee_i) < RTHD \qquad (4.3)$$

If *RTHD* is sufficiently small this ensures the validity of (2.2)…(2.6) throughout the computation. Is this a sufficient condition to also have $ee_i \approx e_i$ for each *i*?

As a measure of accuracy of the estimated error we define the ratio *k* between the true error *e* and the estimated error *ee*:

$$k = e/ee \qquad (4.4)$$

In principle *k* can assume any value, in the absence of some kind of *constraint*, therefore we should find a method capable of rendering very improbable, or almost impossible, the events of large values or of negative values of *k*, in order to have *ee* of some usefulness. Ideally *k* should be 1, but in practice we can tolerate a small spread of values around 1. Our purpose is, in fact, to use the estimated error *ee* to determine a confidence interval for the true error *e* and thus for the true value **x**.

From (4.4) we have: $e = k*ee$. If we knew the statistical distribution of *k* and that, in any case, it holds a condition of the kind $k_{min} \leq k \leq k_{max}$, in the sense that:

$$P(k \in [k_{min}, k_{max}]) = 1 - \alpha \qquad (4.5)$$

with $\alpha \approx 0$, so that it is very unlikely to have values of *k* outside this range, in case $ee \geq 0$, this implies $k_{min}*ee \leq e \leq k_{max}*ee$, which corresponds to a confidence interval for the true value **x**:

$$[x+k_{min}*ee, x+k_{max}*ee] \qquad (4.6)$$

---

[4] Although the derivation of equation (4.2) should be clearly justified at this point, after the analysis of all possible cases made above, one of the reviewers observed that *"the proof of (4.2) is not obvious"*. He/She observed that the case: *"x = -ee and |x| not small"*, would invalidate it. Well, this is something which cannot happen in practice. In fact we should remember that *ee* =0 at the beginning. If during a computation $ee \rightarrow -x$, and $x \gg 0$, this would imply that $|ee/x|$ would exceed *RTHD* well before $|x+ee|/EEZ$ approaches 0. Thus (4.2) remains valid, under the assumed conditions.

The reviewer proposed an alternative definition: $re(x, ee) = |ee|/(|x|+EEZ)$. This definition approximate (4.2), except for case 3), where it tends to be perhaps excessively severe. In fact, if $|ee| \approx |x|$, the relation $|ee|/(|x|+EEZ) < RTHD$, reduces approximately to $|ee|<EPS$.

with a confidence level of 1 - α (of course, if *ee* < 0, the interval should be reversed). Clearly we would like to have α very small, e.g. $\alpha \leq 10^{-5}$, so that the confidence level of the interval (4.6) is almost 1, while keeping this interval as small as possible, i.e.

$$k_{max} - k_{min} = \Delta k \leq M_k \qquad (4.7)$$

where $M_k$ is a small integer. If we can ascertain the validity of condition (4.5) in practice, this would imply that the proposed method is effective. We will discuss the experimental results in the following paragraph, but let us make some considerations first.

Can we evaluate the likelihood of having *k* « 1, or *k* » 1? To be able to express quantitative evaluations we should also have an estimate of the error on the error estimation, which gives rise to an endless process. Another possibility is to interpret *ee* as a stochastic variable and try to estimate its variance. Our attempts in this sense, according to a method outlined in [1] and [23], failed to produce realistic estimates. The variance tends to be overestimated because we cannot take into account the compensation of errors, so that it can only monotonically increase, producing results which are too pessimistic. Errors, as stated above, compensate quite frequently in practice. The reason for this is that the individual errors generated at each step are not always independent; they are frequently correlated instead, because variables and intermediate results are frequently reused in subsequent computations. Thus it also turns out that the distribution of *k* is not normal (as we will see later), because the conditions for validity of the central limit theorem do not hold. However the qualitative observations and the statistical results that we are going to discuss in the next section support relation (4.7) and show that it holds in practice.

It is not rare to come across an exact calculation of *x*. This happens when we can represent exactly the values involved in our computation in the FP system we are using, or when the roundoff errors balance each other exactly in the final result. In this case *e* is 0, if we also have *ee* = 0 we can assume *k* = 1, otherwise, if *ee* ≠ 0, this implies *k* = 0, therefore it must be:

$$k_{max} \geq 1, k_{min} \leq 0 \implies M_k \geq \Delta k \geq 1 \qquad (4.8)$$

We can tolerate *e* = 0, *ee* ≠ 0, as long as *ee* is small, so that the interval (4.6) is small too:

$$e = 0, |ee| < M*EPS \qquad (4.9)$$

where *M* is a small integer and *EPS* is a very small number as defined earlier.

Having *ee* = 0 while *e* ≠ 0, so that *k* = ∞, is fortunately a very improbable event, as it would require exact cancellation (at one or more steps) in the computation of *ee*. In fact, in our statistical observations, it never occurred. Although critical, but arbitrary, ad hoc examples can be provided easily[5], and although one can be worried by cancellation taking place according to Murphy's law, it should be stressed here that the proposed method is effective as long as the statistical behavior of error estimates is satisfactory.

## 5. Statistical behavior of error estimates

The statistical distribution of values of *k* has been investigated in a number of test cases. In the test problems the true exact answers should be known a priori, in order to be able to compare them with the computed values. A suitable test set has been derived from the pentagon in-out iteration problem, first proposed by Dobkin and Silver in [24]. This test is also used in [12], and it is described in [7] and [18] too. The problem consists in finding the intersection points of the diagonals of a convex pentagon, which form an inner pentagon. Extending the sides of the inner pentagon, and computing their intersections, we should get the original pentagon, unless roundoff errors are involved.

The problems of our basic test set have been formed with a parametric pentagon *P* of vertexes: {(0,0),(1,0),(1+*delta*,1),(1,1+*delta*),(0,1)}. An instance of the pentagon *P* for *delta*=0.2, is shown in Fig. 1. The figure also shows the pentagon *in(P),* which is formed by the intersection of the diagonals of *P*. The operation *out* is performed by extending the sides of a pentagon and forming an outer pentagon whose vertexes are at their intersections. In theory it is clearly *out(in(P))* = *P.* The operations *in* and *out* can be iterated a number of times and, as long as they are balanced, they give as a theoretic result the initial pentagon: *out$^n$ (in$^n$ (P))* = *P,* for any integer *n* > 0, which we call the *depth* of the process.

In our test set we considered the operations *out(in(P)), out$^2$ (in$^2$ (P)), out$^3$ (in$^3$ (P)),* in order to form problems with a computation sequence of variable length. We generated a hundred of these problems for each depth of iteration, each time taking a random value of *delta,* respectively in the range: [0, 0.001], [0, 0.01], [0, 0.1].

---

[5] For instance, one of the reviewers, to prove the unreliability of the method, provided this example: $x=1$, $e_x=10^{-6}$, $ee_x=1.001*10^{-6}$; $y=1.8$, $e_y=1.002*10^{-6}$, $ee_y=1.001*10^{-6}$; $z=x-y \Rightarrow ee_z=0$, but $e_z \approx -2*10^{-9}$, which he/she said *"it's by no means small"*. Well, this is relative, as we have discussed in section 4, but, most of all, this is an ad hoc example with arbitrary numbers. Is this situation likely to occur in practice? The reviewer then insisted: *"The probability of this happening is not negligible, especially when the precision of ee is short and that many computations involve numbers of comparable magnitudes"*.

This is an a priori affirmation, not based on facts. The experimental results that we are going to examine will show instead that this probability is negligible in practice, under proper operating conditions. One of these conditions is that an adequate precision is reserved for *ee*. Why should we make it short, favoring critical conditions, when we can avoid them?

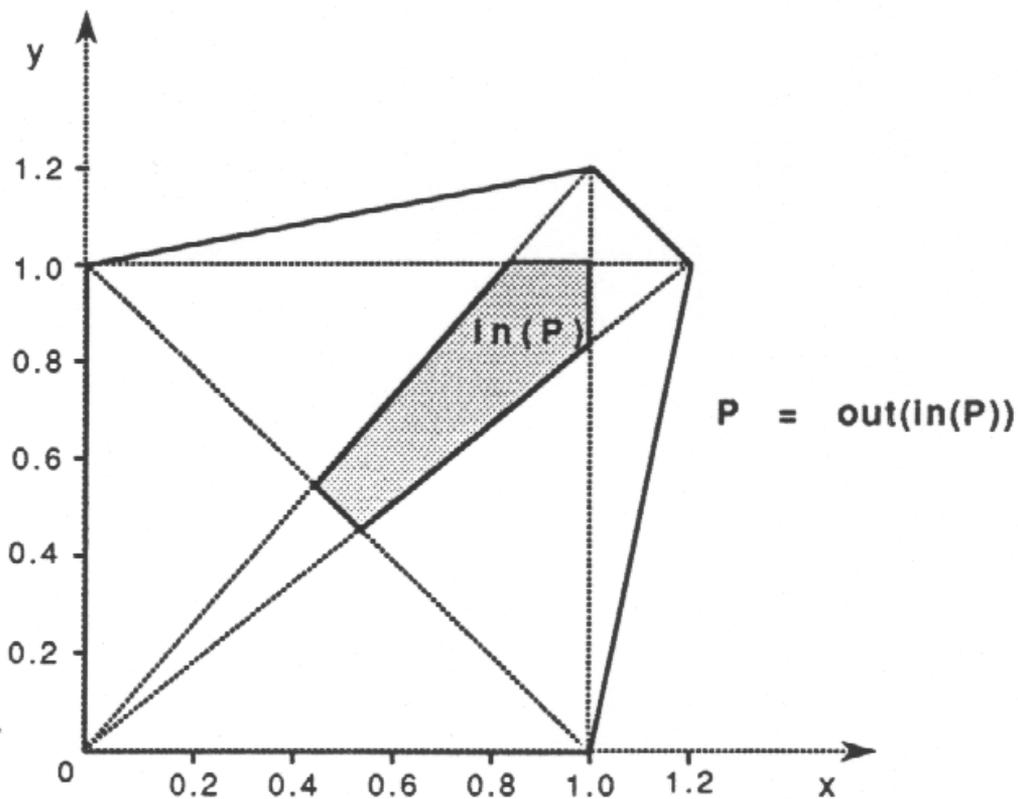

**Fig. 1.** Example pentagons *P, in(P)* for *delta = 0.2*.

The range of variation of *delta* was deliberately adjusted within a limit region, so that, at one end, the larger values of *delta* generate well conditioned problems, while, at the other end, the smaller values generate ill-conditioned problems. Clearly, the deeper the iteration is, the larger *delta* should be, to form a well conditioned problem; this is why we chose larger ranges of variation for increasing depth of iteration. The test trial proceeds repeating the entire sequence, after having translated the original pentagon by the vector [–1,–1], in order to have the "perturbed" vertex close to the origin, and finally by the vector [π, √2], in order to have all the coordinate values larger than 1, and not representable exactly in a FP system. Instances of the test pentagons, for each of the three basic locations, are shown in Fig. 2. The sample that was studied is thus formed by 900 problems, each consisting of the calculation of five coordinate couples for a total of 9000[6] resulting numbers and corresponding values of *k*. In this way we tested the behavior of *k* in a complete span of situations, including very critical conditions.

---

[6] Well, with today's computers, increasing these numbers by at least an order of magnitude would be easy, thus providing more accurate statistics for the results, but we should consider that the original tests were performed on a Macintosh SE! Nevertheless these numbers are high enough to provide meaningful statistics.

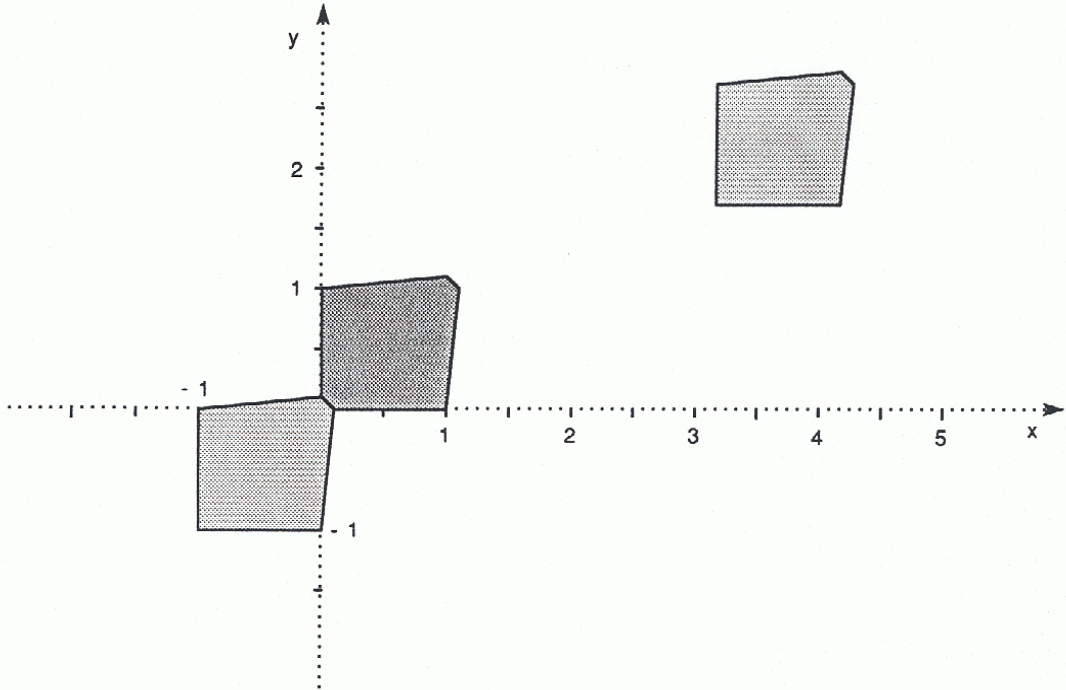

**Fig. 2.** Locations of base pentagons in the test set.

The errors that accumulates in the final results increase with the condition number of a calculation. For large condition numbers, i.e. for ill-conditioned calculations, the errors become very large and tend to assume a random behavior. What happens to a resulting value of a calculation happens to the error estimation as well, so that, at a certain point, also the error estimation becomes very noisy and eventually meaningless. How can we realize this situation? Or, in other words, how can we detect ill-conditioned problems that would lead the error estimation to diverge? This is a main point.

To solve the problem posed by the previous questions, we have to introduce here another mainstay of our method. We have found that an effective approach to control the divergence of error estimation is to monitor the *maximum relative error*, given by (4.2), during a computation. This means that a FP value should be accompanied (although not necessarily in an explicit way, as we will discuss in section 12), in addition to the error estimation, by the maximum relative error $re_m$ which has been registered during its evaluation. We thus define an entity *fpe* represented by the triple:

$$(x, ee, re_m) \qquad (5.1)$$

containing this information. As the process of error accumulation is gradual, we cannot know beforehand if a calculation is ill-conditioned, but if, at a certain point, $re_m$ exceeds an opportune threshold *RTHD*, this would sound as an alarm, telling us that *ee* is probably about to get noisy.

If this case should arise we could either: a) discard the problem reporting a failure, b) signal to the user the poor significance of the results, or c) perform an automatic *dynamic extension* of precision and re-execute the latest operations with greater precision [12]. If this last feature is implemented, we have a tool to guarantee that all calculations are performed with $re_m < RTHD$, so that this is also true in the final results, i.e. we can compute results within given error bounds.

The value of *RTHD* is not critical, but it has to be chosen properly, so that condition (4.7) is likely to hold for calculations where:

$$re_m < RTHD \qquad (5.2)$$

i.e. for calculations where relation (4.3) holds. This means that, in the opposite case, the calculation cannot be performed, with the precision we are using to represent our numbers, without incurring the risk of obtaining gross errors, both in the computed values and in the error estimation. This is confirmed by the experimental results shown in Fig. 3.

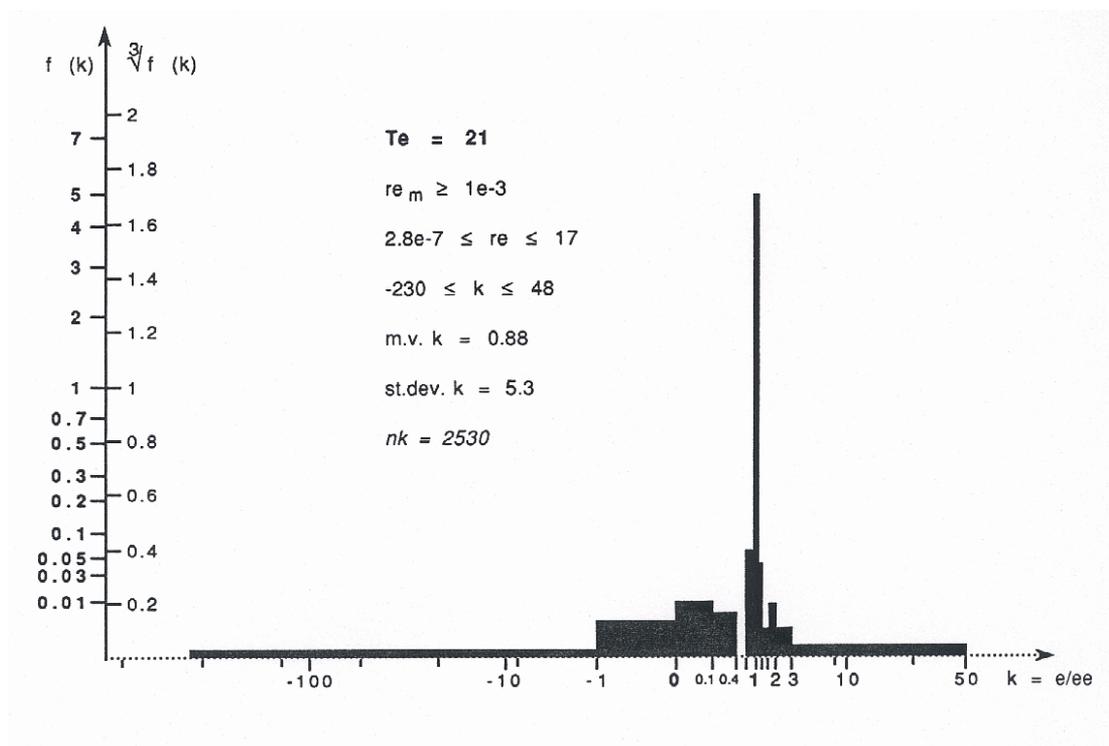

**Fig. 3.** Statistical distribution of *k* for ill conditioned problems, where $re_m \geq 10^{-3}$.

This is the histogram of the observed frequency *f(k)* of the values of *k* for the problems in which, during the computation, we had a value of $re > 10^{-3}$. The mean value *m.v.* and the standard deviation *st.dev.* of *k,* are also given and *nk* is the total number of observed values. Owing to the large domain of *k* we have used a cubic scale for the abscissa. A cubic scale is also used for *f(k)*, to put in evidence the behavior of the distribution in the wide regions of very low frequency density. As shown, the values of *k* are dispersed in a very wide range and their variance is quite high as well. The fact that we have such large values of *k* and of *re*, and negative values of *k* too, signifies that, for some of these problems, both the final values and the error estimations that have been computed are completely noisy and meaningless.

It should be pointed out that, in theory, we could have in some strange case, $re_m < RTHD$, but $k \gg 1$, or $k < 0$. These, however, should be very improbable events, if a proper choice of *RTHD* is made. In fact, what should happen to cause k to diverge? As in a computation the error is initially 0 and then accumulates gradually, the condition (5.2), if *RTHD* is sufficiently small, guarantees that at least the first computations are executed accurately. Thus, if in a sequence of calculations we always have $|ee| < RTHD*|x|$ (which implies $re_m < RTHD$), but at a certain point in the sequence $|e| \gg |ee|$ or $e/ee < 0$ (which implies a divergence of k), this means that the errors on *ee* constantly reinforce themselves, rather than balance each other, at least partially, and also that they are all in the opposite direction to the contribution of the local error *le*, provided at each step, and great enough to inhibit it, in order to have the resultant quantity always confined within a small range. The concomitance of all these circumstances appears statistically quite difficult to occur.

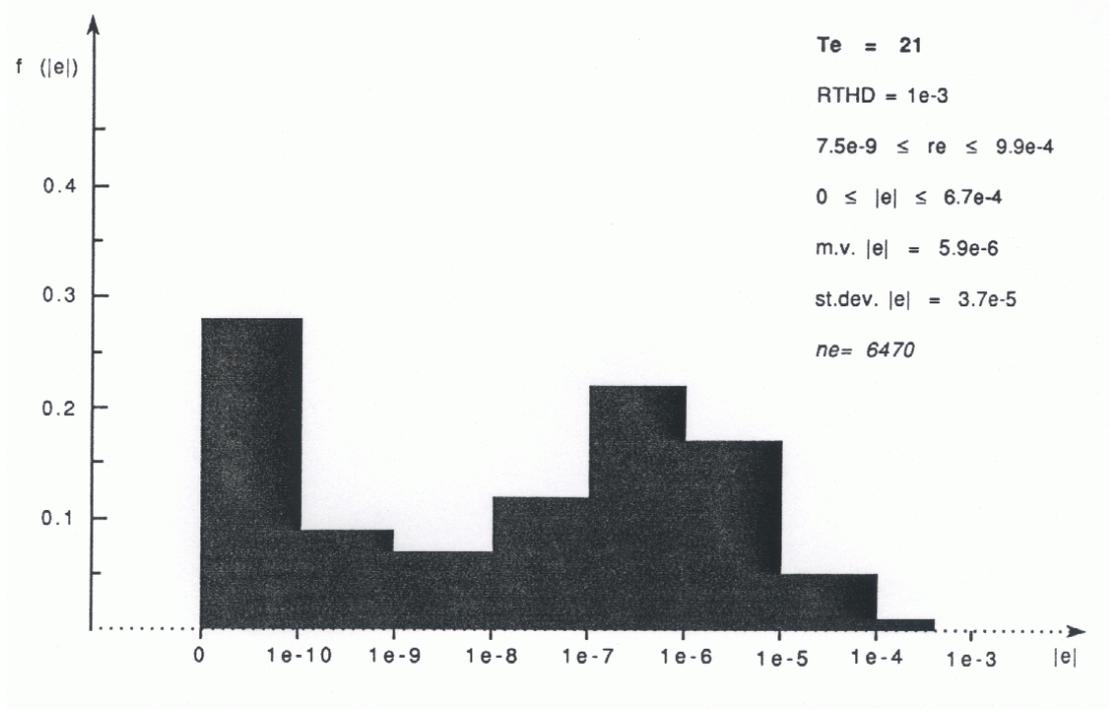

**Fig. 4.** Distribution histogram of the absolute value of the true error |e| for problems with $re < 10^{-3}$.

To verify these arguments and quantify our qualitative observations we are now going to analyze the results of our experimental tests. We fixed five thresholds $RTHD_j$, $j = 1, 2, \ldots, 5$, in the range $10^{-5} \leq RTHD_j \leq 10^{-1}$, and we separately collected the statistics of those problems where $re_m < RTHD_j$. Because of the particular nature of the problem that we used, the final values of the coordinates are strictly correlated, hence if a problem shows a large value of $re_m$ in just one coordinate, this is enough to classify that problem as ill-conditioned. While different values of RTHD were tested, the value of EEZ, was kept fixed. This is equivalent to using different values of EPS for each value of RTHD. We set $EEZ = 10^{-6}$, which corresponds, for instance, to $EPS = 10^{-10}$ when $RTHD = 10^{-4}$. Further on we will discuss the influence of RTHD on the distribution of k, and we will examine in more detail, as a typical situation, the results we had with $RTHD = 10^{-3}$ in various operating conditions. Indicatively, Fig. 4 shows the distribution of the absolute value of the true error |e| with $RTHD = 10^{-3}$, and $T_e = 21$.

The correct choice of the value of RTHD depends primarily on the precision used in the representation of FP numbers, i.e. on the number of bits in the mantissa, and, to a lesser extent, on the particular application and the range of meaningful values involved. In our experimental tests, we have used a FP format with a mantissa of $T = 31$ bits to represent the value. For the representation of the error estimation we have tested different FP formats with a mantissa of $T_e = 31, 21, 16$ bits.

Let us see how the statistical distribution of k depends on the parameters RTHD, and $T_e$.

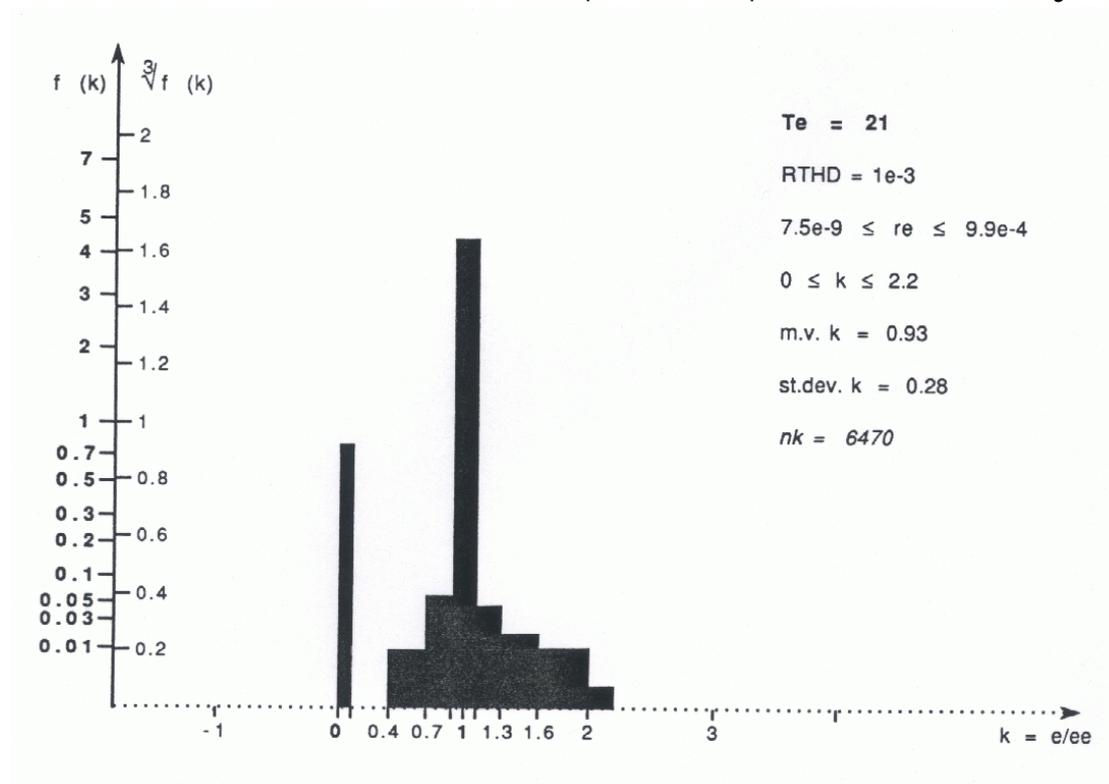

**Fig. 5.** Observed statistical distribution of k in a typical case placing the threshold of re at $10^{-3}$.

Fig. 5 shows the histogram of the observed frequency distribution of values of k for our typical case: $RTHD = 10^{-3}$, and $T_e = 21$. The distribution is bimodal, according to our expectations. Most of the values are densely collected around 1, with some exceptions scattered up to 2.2, but we also have a peak in correspondence of $k = 0$. It is worth noting that we have no instances of *k* < 0, or *k* » 1 here. By comparing this result with the distribution of Fig. 3, where an upper bound on $re_m$ is not enforced, it is evident that *constraining* the computations by $re_m$ < *RTHD* is effective in hampering a divergence of *k*.

The statistical distribution of the values of *k* for $T_e$ = 21 is not particularly sensitive to changes of *RTHD*, at least in the domain $10^{-5}$ ≤ *RTHD* ≤ 0.1, which we considered. The domain of *k* decreases to 0 ≤ *k* ≤ 2 for *RTHD* ≤ $10^{-4}$, but no big changes are observed in the distribution. The variability of *k* is even lower for $T_e$ = 31, due to the increased precision. In this case we always have 0 ≤ *k* ≤ 2 if *RTHD* ≤ 0.1; however, as the threshold increases, the frequency of 0 values decreases, because it is less likely that a value is computed exactly. The dependence of the distribution from *RTHD* is, on the contrary, particularly evident with $T_e$ = 16, which is a symptom of insufficient precision in representing *ee*.

We collected our observations in the graph of Fig. 6, which shows the frequency of *k* as a function of *k* and *RTHD*. The areas with higher frequency densities are darker than areas of lower density. The most important thing to observe in the graph is the dependence, from *RTHD* and $T_e$, of the extreme values of *k*, which define the interval [$k_{min}$ , $k_{max}$]. It appears evident that if we allow for a sufficient precision in the representation of *ee* and we choose a sufficiently small value of *RTHD,* the values of *k* are distributed within a very small range. For *RTHD* ≤ $10^{-4}$ and *Te* ≥ 21, we have 0 ≤ *k* ≤ 2, which is quite satisfactory[7].[8]

---

[7] So, is it clear at this point that it's all a matter of conditioned probability? We can have the true error confined in a small range, around the estimated error, with a very high probability, as long as the relative error is maintained below an appropriate threshold and the estimated error is computed with sufficient precision. Nevertheless one of the reviewers surprisingly observed that the method: *"treats the error as a known quantity instead of an uncertainty. Consequently the method ... becomes equivalent to simulating extended precision".* The reviewer completely overlooked the statistical nature of the method and missed its practical consequences! He/she also missed the point that the purpose of the method is to provide error estimation of the results, not to increase available precision. Increasing precision doesn't provide us error estimation!

[8] One of the reviewers pointed out that: *"The experiments demand knowledge of the exact results. How can kmin, kmax be estimated in general?".* This is a good question. It is answered in section 7.

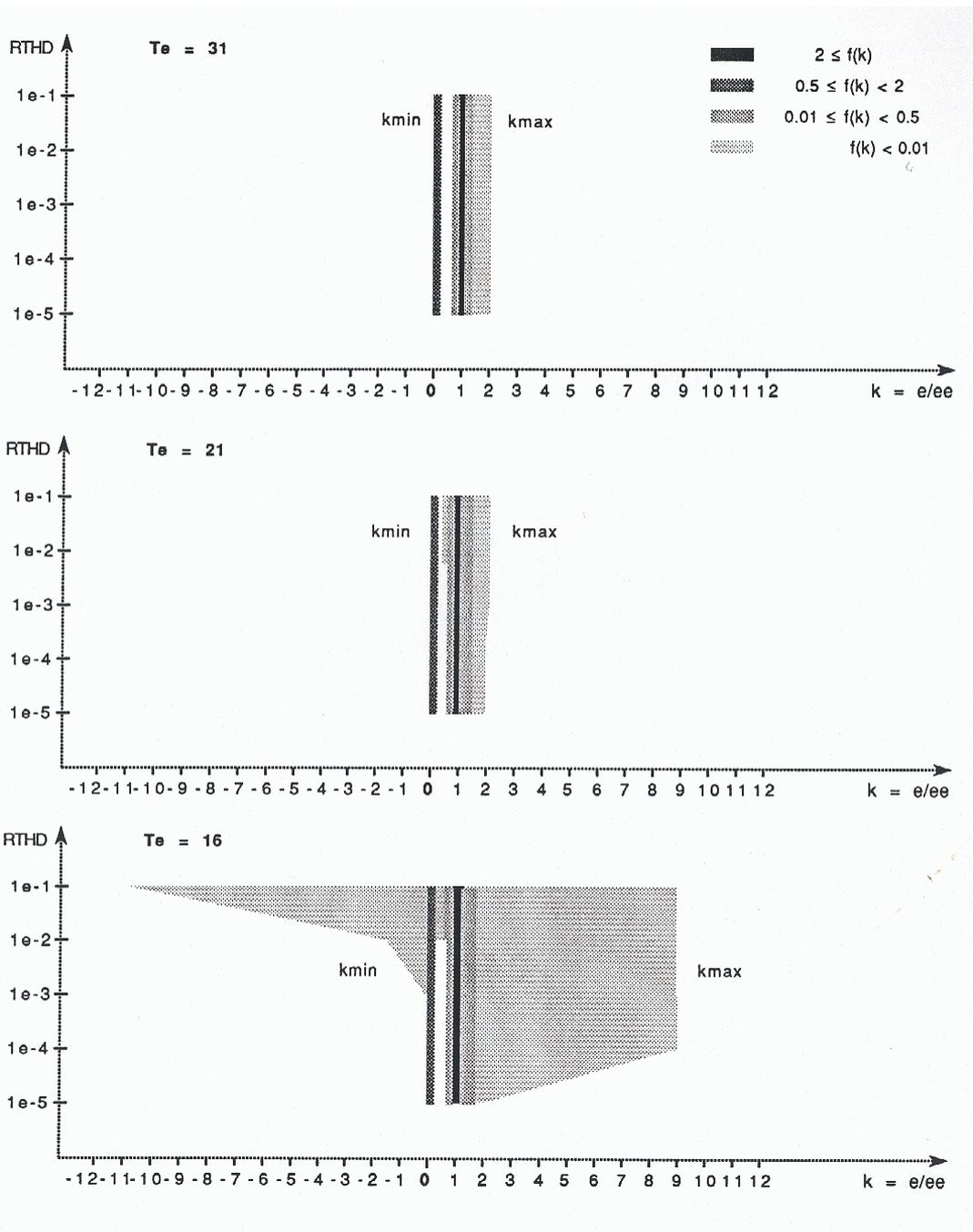

**Fig. 6.** Dependence of *f(k)* from *RTHD* and *Te*.

We collected other data relevant to our reasoning. First we have to verify the statistical validity of (4.9). This is confirmed by the distribution of |*ee*| under the condition *e* = 0, which is shown in Fig. 7.

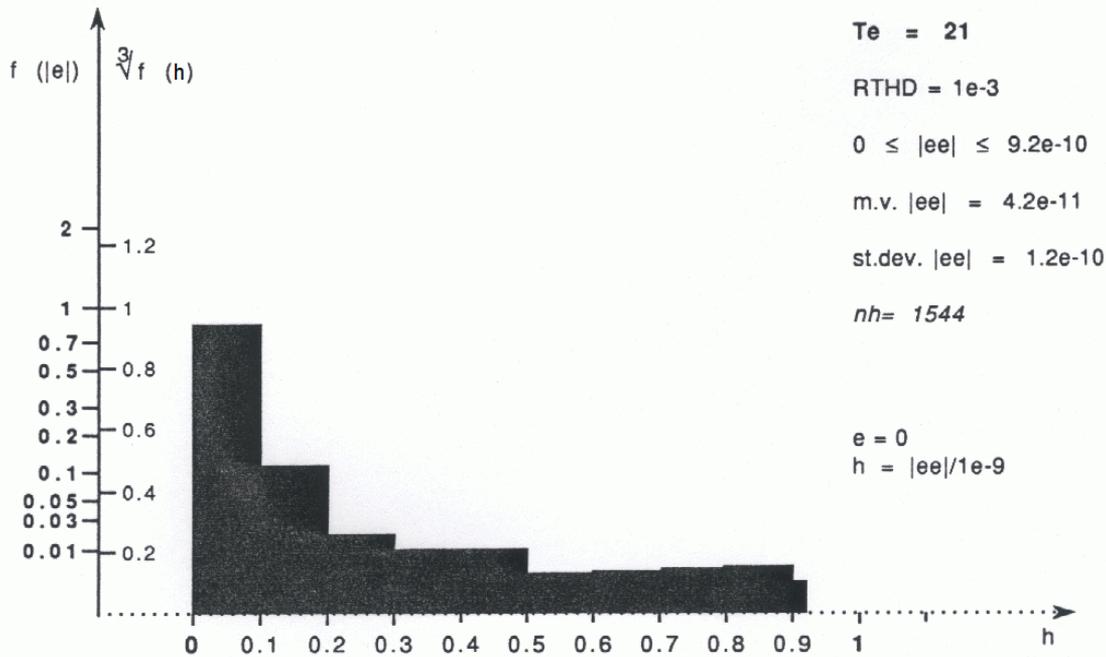

**Fig. 7.** Distribution of the estimated error *ee* when the true error *e* = 0.

It is reassuring to see that the case *e* = 0, *ee* = 0, is very frequent in practice. Also, given *EPS* = $10^{-9}$, which corresponds to the case *RTHD* = $10^{-3}$ at hand, the average value of |*ee*| is about 0.04\**EPS*, and in any case we have |*ee*| < *EPS*, which is more than satisfactory. The fact that sometimes |*ee*| approaches *EPS* tells us that we have not chosen a value of *EPS* which is too big. The distribution of |*ee*| does not change much, when *RTHD* is varied and $T_e \geq 21$. With $T_e$ = 16, instead, we have much more noise in *ee*, the variance increases significantly and the variability of |*ee*| extends up to 2.4\**EPS*.

In section 4 we also made the supposition that |*ee*| < *EPS*, in the case *x* = 0, *x* = 0. We observed |*ee*| < 0.0012\**EPS*, in the case corresponding to Fig. 7, which means that the chosen value of *EPS* safely accommodates for the largest observed values of |*ee*|, when *x* = 0, *e* = 0, *ee* ≠ 0. This result is almost insensitive to variations of our parameters. This can be explained by the fact that in this case the computations are quite accurate.

The actual value of *EPS*, i.e. a correct choice for it, depends on the precision used in the representation of FP numbers, and on the range of values involved in the applications. In the test sample we have a range of significant values about five orders of magnitude below and one order of magnitude above unity. In these operating conditions, and with the FP format we used, the choice which we made, based on the relation *EPS* = *RTHD*\*$10^{-6}$, appears to be quite satisfactory.

January 2012					22

## 6. Variation on a theme

In which circumstances do we record the larger values of *k* when $re_m$ < *RTHD*? A closer look at the results shows us that, the larger deviations of *ee* from the true error *e* are found for some of the smallest values of *e*. This is probably due to the loss of precision that occurs when the most significant "extra bits" in the accumulator, i.e. in the *v* part of (3.4), are set to 0. The fact that in these circumstances *e* and *ee* are small, suggests that we study the distribution of the variable *c* given by:

$$c = e/(sign(ee)*(|ee|+QEPS)) \qquad (6.1)$$

where *QEPS* = *Q*EPS,* and *Q* is an opportune factor. The definition (6.1) of *c* differs from the definition (4.4) of *k,* in that the simple ratio *e/ee* is "perturbed" by the amount *QEPS,* thus c is defined also for *ee* = 0, and the modulus of *c* is bounded to |*e/QEPS*| when *ee* → 0. This has the consequence however of setting to *QEPS* the width of the smallest confidence interval of the error. In fact we have: *e* = *c*sign(ee)*(|ee|+QEPS)*, so that, setting *ce* = *sign(ee)*(|ee|+QEPS)*, if *ce* > 0, we have: $c_{min}$**ce* ≤ *e* ≤ $c_{max}$**ce* , which corresponds to a confidence interval for the true value ***x**:*

$$[x+c_{min}*ce , x+c_{max}*ce] \qquad (6.2)$$

Recalling (4.8), it is $c_{min}$ ≤ 0 and, in this case too, it is likely to have $c_{max}$ ≥ 1, so it follows that the interval (6.2), is necessarily greater than, or equal to [*x, x+QEPS*] (or [*x-QEPS, x*] if *ce* < 0). Therefore the value of *Q* is chosen as a compromise between the demand of having *QEPS* as small as possible and that of limiting the value of $c_{max}$ . In our experimental test we set *QEPS* = $3*10^{-10}$, which implies $Q = 3*10^{-4}$/*RTHD*.

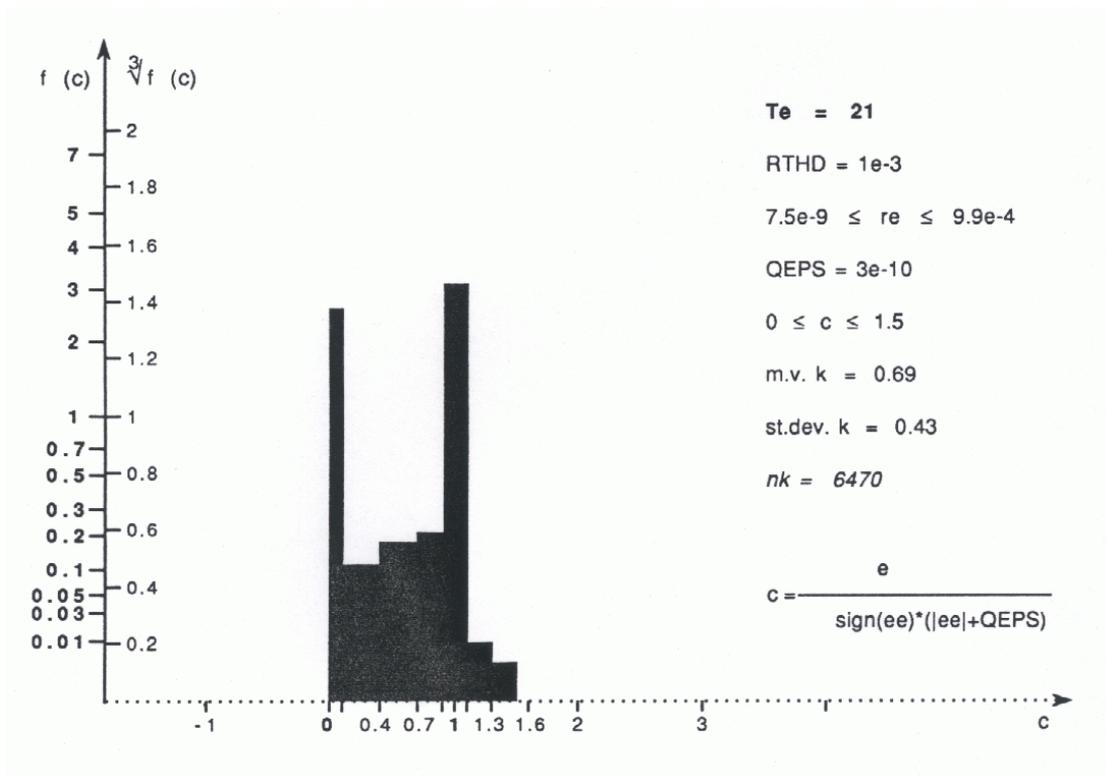

**Fig. 8.** Distribution of c in our typical case.

Fig. 8 shows the histogram of the distribution of $c$ in our typical case. It can be seen, as expected, that $c_{max} - c_{min} = \Delta c < \Delta k$. We have in fact $0 \leq c \leq 1.5$. The values of $c$ are spread over a smaller interval with respect to $k$, and they are more uniformly distributed over this interval. In fact we can observe a considerable increase of frequency around 0, as this time $c = 0$ catches all the cases in which $e = 0$, regardless of the value of $ee$. The gap existing in the distribution of $k$ between 0 and 1, is filled up in the distribution of $c$ with values whose probability increases approaching 1, then the probability of having $c > 1$, drops drastically.

Fig. 9 sums up the dependence of the observed frequency $f(c)$ on $RTHD$ and $T_e$. Comparing it with Fig. 6, we can notice the general reduction of the interval $[c_{min}, c_{max}]$, with respect to $[k_{min}, k_{max}]$, even for the case $T_e = 16$.

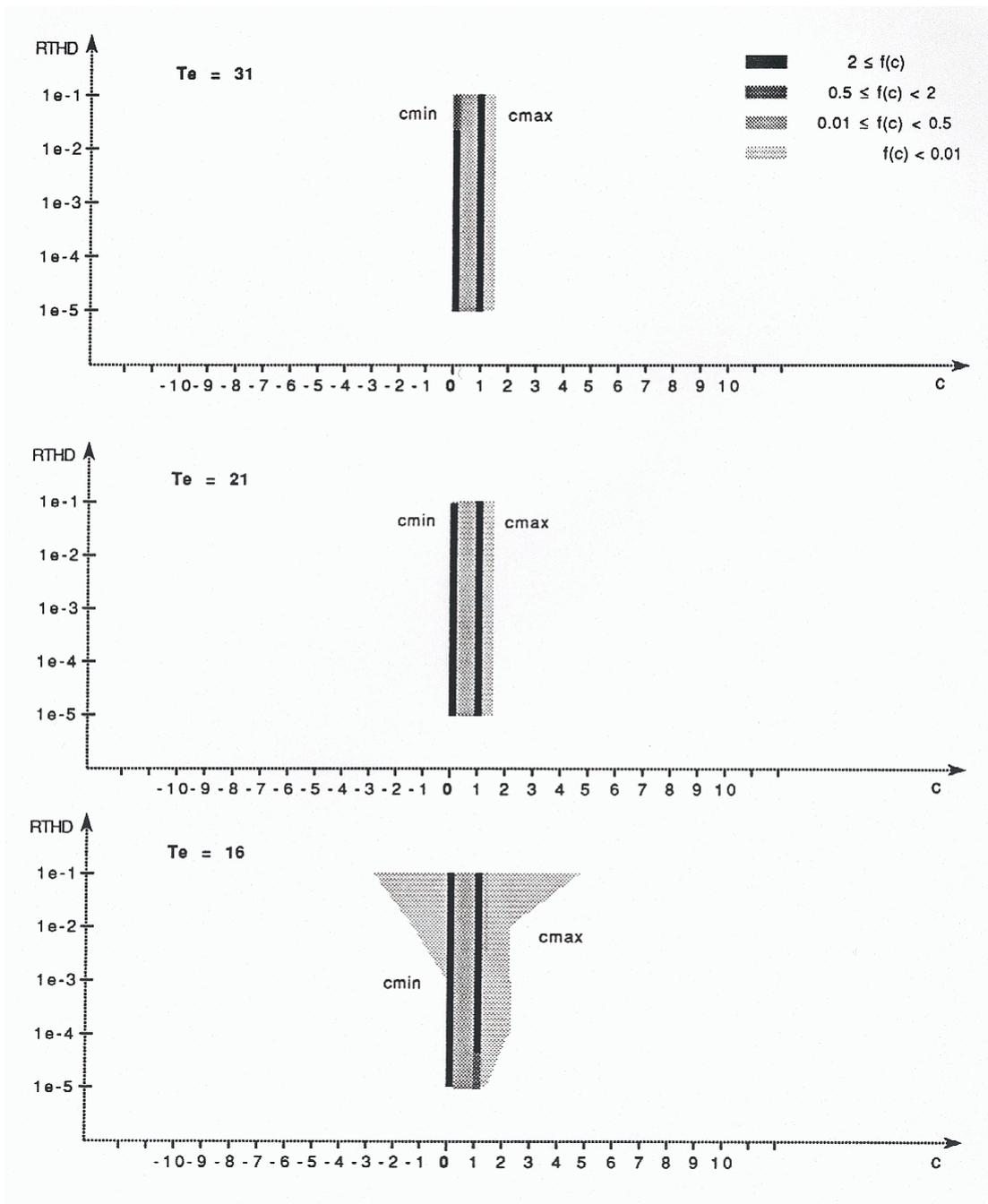

**Fig. 9.** Dependence of *f(c)* from *RTHD* and *Te*.

Justification of the opportunity to use *c* rather than *k* depends on the application: *c* has a more compact distribution with respect to *k*, but it has some drawbacks in dealing with small values of *x*, since the intervals of confidence, associated with these values, will be widened by the influence of QEPS, like observed above.

## 7. Confidence intervals for *k* and *c*

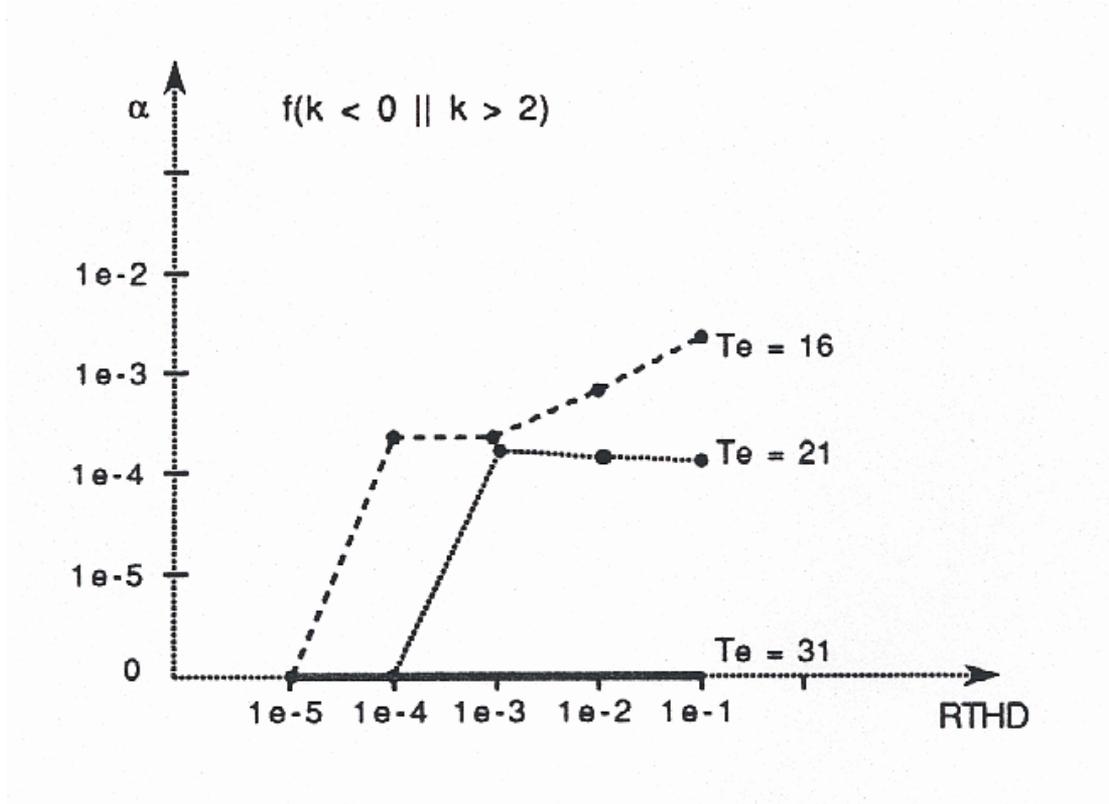

**Fig. 10.** Observed frequency of *k* outside the interval [0, 2].

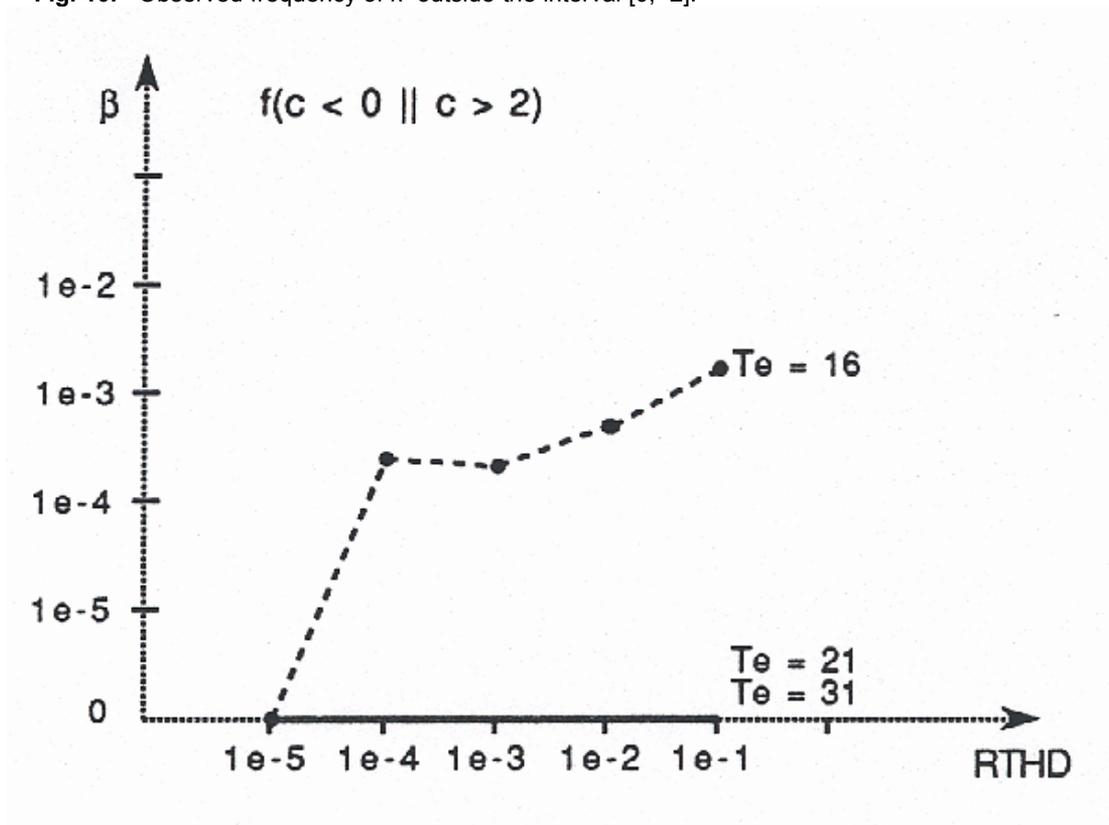

**Fig. 11.** Observed frequency of *c* outside the interval [0, 2].

The graphs of Fig. 10 and Fig. 11 summarize and compare the estimated probabilities $\alpha$ and $\beta$ of having respectively *k* and *c* outside the interval [0, 2]:

$$\alpha = P(k \notin [0,2])$$
$$\beta = P(c \notin [0,2]) \tag{7.1}$$

The graphs are built from the observed frequencies of the two events and display also the dependence of $\alpha$ and $\beta$ on $T_e$ and *RTHD*. The interval [0, 2] appears as a reasonable assumption for [$k_{min}$, $k_{max}$] or [$c_{min}$, $c_{max}$], based on the observed statistic for *k* and *c* (although a smaller interval, like [0, 1.6], could have been chosen for $\beta$).

The use of the method, proposed in this document, can be effective only if the probabilities (7.1) are very low. Ideally they should be 0, but values below $10^{-5}$ should be acceptable for most applications. The estimated confidence level $C_k$ of having $0 \le k \le 2$, i.e. [$k_{min}$, $k_{max}$] ≡ [0, 2], is thus $C_k$ = 1 - $\alpha$, and analogously we have $C_c$ = 1 - $\beta$ for the variable *c*. In our simulation, for $T_e$ = 21 and *RTHD* = $10^{-3}$, we have $C_k$ = 1 - $3.7*10^{-3}$ = 0.9963, and $C_c$ = 1. Such an high value for $C_c$ in obviously suspicious, most likely a greater sample, than the one we used in our test trials, is required to estimate $C_c$ more precisely.

We carried out several other test trials, conceived to be more severe and exhaustive than the basic test described earlier, with the purpose of ascertaining the validity of the, quite positive, results of above. The in-out iteration problem was performed on thousands of randomly generated pentagons, and with different depth of iteration. In some test trials the iterations were also cycled in order to build up sequences of variable and longer length. The results were in substantial agreement with those of the basic test.

The statistical distribution of *k* seems to be little influenced by the length of the sequence of computations. It is greatly influenced instead by the distribution of the condition number of the problems in the trial. In fact, if the problems of the trial are mainly well conditioned, the probabilities $\alpha$, $\beta \to 0$. These probabilities increase with the frequency of ill-conditioned problems in the sample. Nevertheless, in all circumstances, we observed that they can be reduced to 0 by increasing $T_e$ and/or decreasing *RTHD*.

Of course the experiments have been carried out in a situation where the exact results are known, in order to be able to collect statistics on the errors. The test set, as stated in section 5, is made of a complete span of situations which may arise in practice, ranging from well-conditioned to ill-conditioned problems. It is intended to be representative of practical situations, even the most critical ones; for this purpose the frequency of ill-conditioned problems is higher than normal. Therefore we should expect a similar or, most likely, a more favorable behavior in general. In fact, we carried out other experiments, with matrices and computer graphics algorithms, which results have not been systematically collected, but in all these cases we have often found more compact distribution of the errors, or at worst in agreement with the statistics presented above.

In conclusion, if the values of α and β are estimated by representative and worse than normal test cases, it appears legitimate and safe to extrapolate the results of the test cases to applications in general. Thus, if α and β in the test cases are very low, $C_k$ and $C_c$ are very close to 1, and we should expect similar or better results in general[9].

On the basis of these considerations we can make appropriate choices for $T_e$ and *RTHD* in a practical system, in order to make acceptable the confidence levels $C_k$ and $C_c$. Increasing $T_e$ and decreasing *RTHD* will increase $C_k$ and $C_c$, but obviously it will cost more in terms of required space and execution time (in case automatic, dynamic extension of precision is enforced, see section 5).

## 8. Space required for *fpe* structures

In an *fpe* structure, as defined in (5.1), the error estimation *ee* and the maximum relative error $re_m$ should always be associated with the value *x* of a FP number. There is thus the demand to limit the space taken by the two error fields. It is evident that the representation of $re_m$ does not need to be accurate, as we can store just a superior limit of the true $re_m$, updating it when it is exceeded, so that we need just a few bits for the mantissa of $re_m$.

Moreover, although very useful in this experimental phase, in practical applications we can avoid storing $re_m$ with each number, as we will see in section 12, so this space can be completely saved. In this case the *fpe* structure reduces to the couple (2.1), in place of the triple (5.1).

Differently, the effects of accuracy in the representation of *ee* are very important, as pointed out in the previous sections. However, we have seen that, as long as we require a small relative error, i.e. if $RTHD \leq 10^{-4}$, varying $T_e$ in the range $21 \leq T_e \leq 31$ does not greatly affect the accuracy of *ee* and consequently the distribution of *k* or *c*. Instead, we have significant variations going down to $T_e = 16$ bits, because, at this point, the loss in accuracy is excessive.

On the basis of these considerations, a possible format to represent *fpe* as a triple is shown in Fig. 12. In this particular solution, which has been adopted in our test trials, 1/2 of the total available space is reserved to the value part, 1/3 to the error part and 1/6 to the maximum relative error part. Obviously other partitions of the available space, which respect the requirements stated previously, are possible.

But, how much space do we globally need in practice? Should we deduce from the above, that when using *fpe* structures we need approximately twice as much space as we need to operate with simple FP numbers? This can be the case, if we reserve for the value part of the structure the same space we were used to assign to FP numbers in our applications, and about the same space for the error part. However, by taking full advantage of the information contained in *fpe* structures, we should be able to implement most applications using a smaller space, i.e. just the amount of space needed to achieve the necessary accuracy in the final results. This should be possible in particular if we use the error information to control a representation of numbers with variable precision, as advocated above. This may entail the use of smaller formats in general and of larger formats only when required.

---

[9] The estimated confidence levels $C_k$ and $C_c$ could perhaps be lowered by a safety coefficient, to take into account hypothetical worse cases that may arise in practice, but this appears as an arbitrary procedure, thus we preferred to avoid it.

**9. Software simulation**

A software simulation can, to some extent, be implemented. Things are easier in an object oriented language like, for instance C++. In this case we can represent objects of type *fpe,* defining the class:

```
class fpe {
   float value;
   float error;
   float relerr;
};
```

Then we provide the constructors for initialization and the methods for doing the arithmetic operations with objects of this type, and for applying functions to them. These methods could be applied to *fpe* objects, in the same way as we do for the types `double` or `float`, by overloading the usual operators and functions.

The values of *RTHD, EPS* and eventually of *QEPS,* could then be adjusted for the particular FP system available in our machine and for the particular application. A signal could be generated if, for one of these objects `x,` we have `x.relerr` > *RTHD,* so that special processing could take place.

The main obstacle to a software simulation, using a high level language, is that we don't have direct access to internal machine registers and to details of arithmetic operations. Thus we cannot monitor the rounding process, unless we exploit the use of type `double,` to know about "extra bits" with respect to the type `float,` which we have used to represent the value in the class `fpe.`

For instance let us suppose we have defined `fpe` as:

```
class fped {
   double   value;
   float error;
   float relerr;
};
```

if `x` is of type `fped` and we want to initialize this object to the value 0.6 (which exact representation is the binary periodic fraction 0.1001 1001… ) we would try:

```
x.value = 0.6;
x.error = x.value - 0.6;
```

but we will not be able to initialize the error field to a value different from 0. In fact constants are internally represented as `double` entities, so they are rounded at the time they are created by the compiler. We then lose track of the amount of this rounding, which is exactly what we need to initialize the error field. The same considerations hold for rounding in intermediate computations. To be able to do anything, we need to represent numbers in `fpe` with less precision than is used to represent constants or to do intermediate computations, this is the reason why the definition of `fpe` works, while the definition of `fped` does not.

Our software simulation has been implemented in C for pure experimental purposes, with no care paid to optimizing time or space. We have used the definition:

```
struct fpe {
   double   value;
   double   error;
   double   relerr;
};
```

In our machine the type `double` has a mantissa 63 bits long. We then wrote routines to manipulate the mantissa of a double internally, thereby simulating a rounding operation to any number of bits in the range $2 \leq T \leq 61$. In our simulation the result of a computation is temporarily available in a `double`, which acts as the accumulator $w$ of section 3. In these operating conditions we have simulated a processor with $T_W = 63$, $T = 31$.

In any case, however, a software implementation will be inevitably slow, compared to simple FP processing, so we believe that the availability of specialized hardware is the only viable solution for the method to be used in practice.

## 10. Using estimated error in algorithms

Having got so far, how should we use these error estimates in our programs? The time for taking advantage of error estimates comes when we have to make decisions based on the computed values. For example, let us suppose that at a certain point of a program we have to make a test of this kind:

```
if(a == b) {
    ...
    }
else ...
```

Recalling the concepts exposed in the introduction, and also considering our experimental results on error statistics, if `a` and `b` are FP numbers this test is likely to fail almost every time, even when the two numbers should correspond to the same real value (to succeed there should be the simultaneous occurrence of exact calculations of `a` and `b`). The decision to be taken in the equality test however may be critical, as singularities play an important role in practice. For example, let `a` and `b` be the semiaxes of an ellipse; if they are equal the ellipse becomes a circle, thus we can take advantage of the special geometric properties of a circle.

The best way of coping with this problem has been, hitherto, to rewrite the test (according to (1.1)) in this way:

```
if( abs(a - b) < EPS*(1+abs(a)) ) {
    ...
    }
else ...
```

The absolute equality has been substituted with the comparison against a relative tolerance. The constant `EPS` is set arbitrarily by the programmer with the consequences that we have discussed in section 1.

Instead, with the proposed approach to numerical computations, using the error estimates, we could be much more accurate. First, we build the two confidence intervals given by (6.2) (or (4.6)) for `a` and `b`, then we check if they *overlap.* If they are disjointed (or touch) there is virtually no possibility that `a` and `b` represent the same number, but if they overlap, it is likely that they represent the same number, so we are allowed to consider them as equal. We could be even more accurate, and consider "how much" the intervals overlap. Using the estimated distribution of probability for *c* (or *k* ), we can calculate the joint probability that the two objects represent or not the same real number, and take a decision consequently; but, obviously, this has a cost. Let us assume that we have a (possibly hardwired) method called `equal` which applies to an object of type `fpe` and takes as argument another `fpe` object. The method builds the confidence intervals (6.2) for the two numbers and checks whether or not the two intervals overlap, so that the two numbers can be considered as equal. Our code may then be rewritten —using a C++ notation— as:

```
if(a.equal(b)) {
    ...
    }
else ...
```

with much more likelihood of correctness.

Let us consider another example: we want to compute the intersection of two straight lines in the plane:

$$a*x + b*y + c = 0$$
$$e*x + f*y + g = 0$$

When should we have to consider the two straight lines as parallel? In theory, we know that we should check whether the determinant $D$ of the matrix of the coefficients of $x$ and $y$ is 0. In practice, it will seldom be 0. Up to now, we have been used to test for $D < \varepsilon$, but what is a reasonable $\varepsilon$? With error estimation, instead, we could test whether 0 belongs to the confidence interval of $D$ (given by (4.6) or (6.2)) and be quite certain that all of, and only, the singular cases will be caught[10].

## 11. Hardware implementation

A special FP processor, implementing the proposed method, will operate on a FP format of the kind shown in Fig. 12, or, most likely, on a simplified structure without the $re_m$ field, as we will see in section 12.

In designing the specialized FP processor a solution has to be found, exploiting parallelism and pipelining, which involves a minimum delay with respect to standard FP processing. To this purpose the equations (2.2)…(2.8) can be rewritten for the total error $te = e_z = pe+le$, pointing out the result $z$ for each operation; we have:

$$te(z = x+y) = e_x + e_y + le$$
$$te(z = x-y) = e_x - e_y + le$$
$$te(z = x*y) \approx y*e_x + x*e_y + le$$
$$te(z = x/y) \approx ((e_x/y) + le) - (e_y/y)*z$$
$$te(z = \sqrt{x}) \approx (e_x/2)/z + le$$

---

[10] The proposed method thus qualifies for absolutely general employment, determining a new paradigm to treat numerical inaccuracies in any kind of algorithm, but geometric problems in particular could benefit from it. It is thus surprising that one of the reviewer affirmed: *"this idea is not powerful enough to solve the robustness issue in a lot of important geometric problems such as solid modeling, hidden line and surface removal etc."*. Exactly the kind of problems for which the method was originally intended for! If I could have used the method (having the necessary hardware support to overcome the performance barrier), I certainly would have taken great advantage of it in my code, realizing much more reliable algorithms than those that can be written using the standard approach discussed in section 1. However the reviewer gave no evidence to support his/her assertion.

It is worth recalling that the evaluation of $le = w-z$ can be made concurrently to other operations as soon as the result $z$ is available! For each case we can make the following observations:

- +, –  In the case of addition and subtraction the evaluation of $e_x+e_y$ can be performed concurrently to the evaluation of $z$. The total error $te$ is provided with one delay to add $le$.
- \*  In multiplication the two parts $y*e_x$ and $x*e_y$ can be evaluated concurrently to $z$, but we have one delay to add them successively, and one delay to add $le$.
- /  In the case of division it is possible to evaluate $z$, $e_x/y$ and $e_y/y$ concurrently, but then we have one delay to multiply $(e_y/y)$ and $z$. The addition $(e_x/y)+le$ can be performed concurrently at this time, then we have another delay to subtract the two resulting terms.
- √  For the square root, the evaluation of $e_x/2$ and $z$ can proceed concurrently, then we have two delays, to divide the first term with the other and to add $le$.

January 2012                                                                                                           33

These observations have been taken into account in the block diagrams of Figs. 13…17, which illustrate the data flow in the various cases.

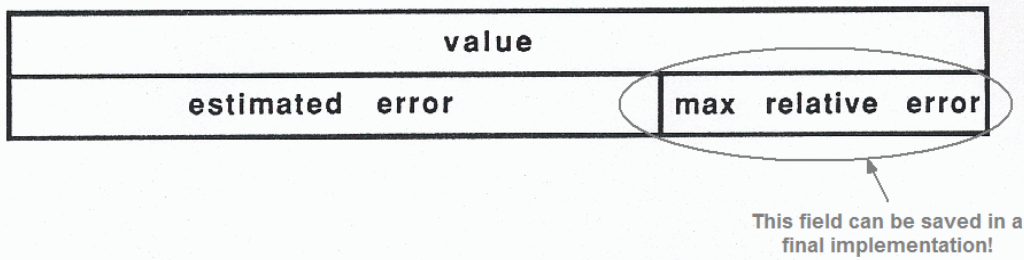

Fig. 12. Format of structure *fpe*.

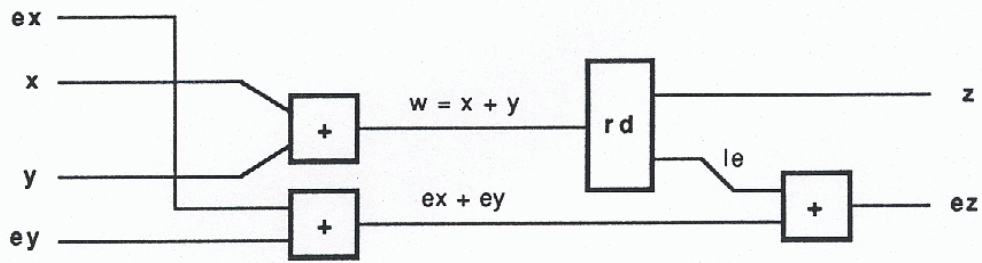

Fig. 13. Exploiting parallelism in addition.

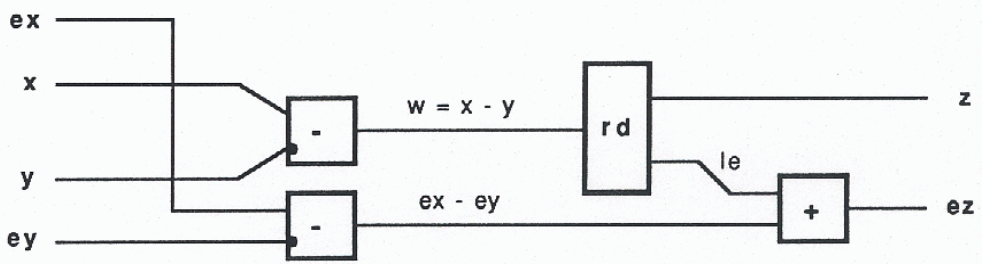

Fig. 14. Block diagram for subtraction.

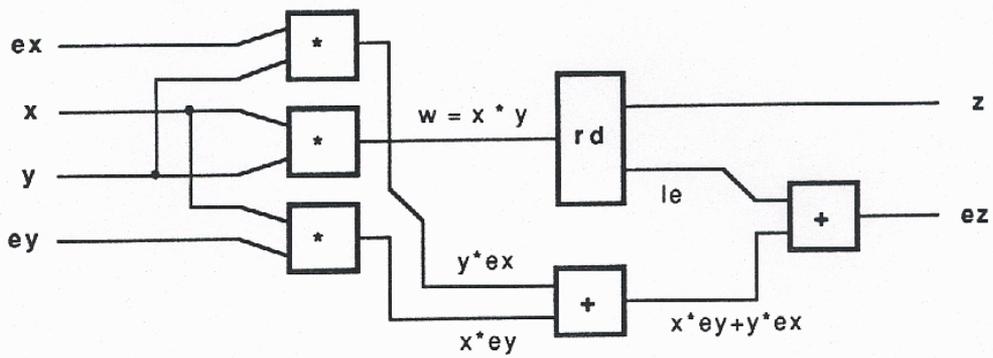

**Fig. 15.** Exploiting parallelism when multiplying *fpe* structures.

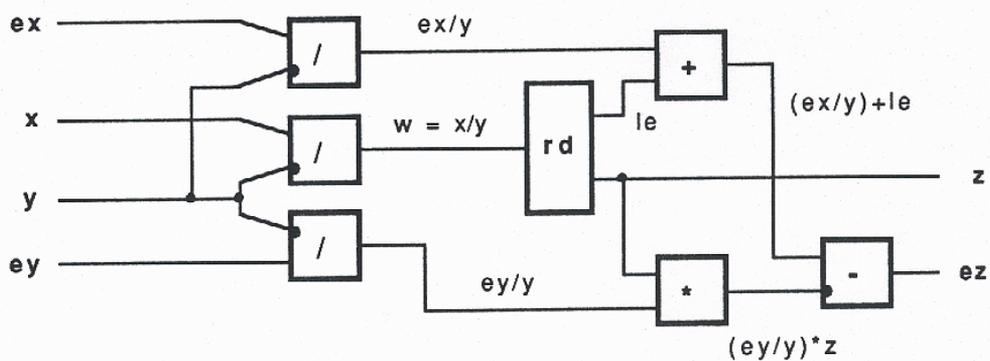

**Fig. 16.** Block diagram for the case of division.

After having computed the value and error estimation parts of the result, we can proceed to the evaluation of the relative error part. This procedure is common to all four arithmetic operations and is shown in Fig. 18.

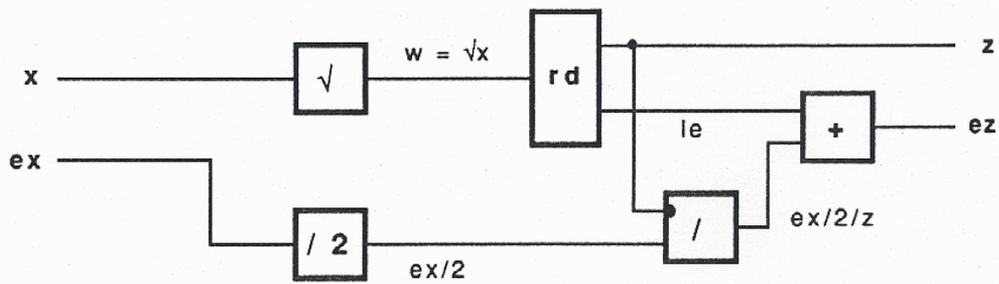

Fig. 17. Block diagram for the square root case.

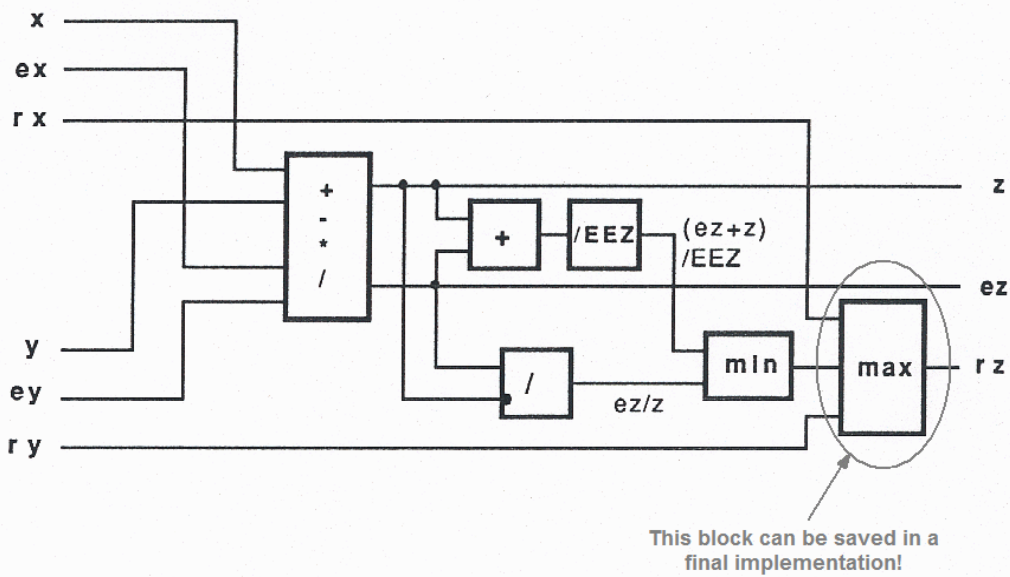

Fig. 18. Complete diagram for arithmetic operations on two *fpe* structures.

The diagrams of Figs. 13...18 can serve as a basis for the design of an arithmetic processor, supporting *fpe* structures, and exploiting all possible concurrencies. A processor of this kind could evaluate arithmetic operations on objects of type `fpe`, with a time penalty factor *TPF* ranging from 3 to 4, i.e. taking approximately from 3 to 4 times longer than standard FP operations. This is far from optimal, because, after performing their task, the modules are left idle till the next operation. The value of *TPF*, stated above, can be deduced from a quick exam of the serial operations that are involved. An exact evaluation will depend on the practical implementation and on some simplifications which can be introduced. For example, examining the diagram of Fig. 18, we can observe that we are not interested in an exact evaluation of $e_z/z$ but only in the most significant digits of the result. Therefore, in a tweaked implementation, the division which appears in the diagram can be faster than a regular FP division. In addition, *EEZ* can be chosen as a power of 2 so that the division by *EEZ* can involve, in practice, only simple exponent arithmetic. But further and greater improvements could probably be achieved by pipelining the operations, in order to maximize the employment of every module.

Finally, this work should be integrated with research in variable precision arithmetic [12], [21], [22], [35], in order to implement *dynamic precision extension*.

The realization of a specialized processor incorporating the method in the hardware, as outlined above, with current technology, should not be anymore a problem and would make the practical adoption of the method feasible for most applications.

## 12. Effects of our proposal and subsequent related work

What happened after the publication of [25]? We must repeat that our work was quite neglected afterward, perhaps also because a CAD review is not exactly in evidence for experts in numerical problems or FP architectures. Nevertheless the paper got some attention around the world [26], [27], [30], [31], [32], [33], [34], although sometimes only for minor aspects of the work. We are thus going to revise and comment only the most relevant publications, we are aware of, that made reference to our original paper.

The most positive appraisal was perhaps that of John Keyser[11] [26], [33]:

> "... When worst-case error estimates are kept (such as in interval arithmetic), one is guaranteed that the final result is located within a certain bounded interval. However, after several computations, this interval may be so large that it is worthless. Masotti's error estimate is much smaller than the worst-case error, and thus the final answer is more useful. In addition, this approach can detect certain ill-conditioned computations. Masotti demonstrates that this approach can be used to perform operations far more robustly than with standard floating point, although the results are still not guaranteed. Masotti's approach also has potential in that it can potentially be implemented in a hardware floating point processor."

Perfect! We only have to point out that, with an appropriate implementation, the confidence level of the error estimate can be made arbitrarily high, so that results are virtually guaranteed.

---

[11] But unfortunately he misspelled my surname in the reference list. Not a big deal, but I guess this can make difficult to find, recognize or count the citation, especially for automatic programs.

The original paper [25] is also taken into detailed consideration in the doctorate thesis of Denise Pirus [27]. Translating from the original text in French, Pirus summarizes our work this way:

> "G. Masotti proposes a method to obtain an accurate estimate of the absolute errors that affect the results. The method consists in calculating, throughout the algorithm, the absolute error of each real. In order to do this, he represents a real number *x* by its value, the estimate of the error associated with it and the maximum relative error *re*. As the error accumulation is gradual, it is possible that, at some point of the algorithm, *re* exceeds a threshold of acceptability, which is dependent on the precision used to represent floating point numbers, if this happens, the user can then be alerted.
>
> The real is thus represented by the triple *(x, ee, re)*. He determines the contribution of the error in a calculation, i.e. the propagated error *pe*, using a Taylor expansion.
>
> To this error he adds the rounding error done by the computer. To determine it he assumes that, at some point, he has at his disposal the result *w* with a number of significant bits greater than the format of the real, and in this case one has:
>
> $z = rd(w)$, rounding function that operates on *w* to obtain *z* in the format stated;
>
> $le = z - w$, local rounding error.
>
> The result of a calculation is thus the triple $(z, e_z, re_z)$, where $e_z = rd(pe + le)$ and $re_z = |e_z|/z$."

This synthesis is pretty much correct, a part some simplification and imprecision in the definition of the various entities, for which it's better to make reference to sections 2, 3 and 4 above. Pirus then adds:

> " The objective of G. Masotti is to use this error estimate to determine a confidence interval for the true error, and thus the true value of the considered real. This is what allows him to compare two numbers more precisely: if *a* and *b* are two numbers to be compared, he builds both confidence intervals and checks if they overlap. If they are disjoint, then *a* and *b* are not the same number, otherwise they are considered as equal.
>
> This method makes it possible to control accumulated errors and to take them into account when comparing real. However, it introduces some problems:
>
> The memory space required to represent a real increase significantly, since, in addition to the number, he stores the absolute error and the relative error."

Pirus seems to have forgotten that not the relative error is stored, but the maximum relative error encountered in the computation ($re_m$ in (5.1)), nevertheless she then comes out with a surprising and very interesting observation:

> "It does not seem interesting to keep this value, because it does not intervene in the following calculations, it is sufficient to compare the relative error with the acceptability threshold, just after the calculation."

This in fact can be true in a final implementation of the method, although, in the experimental phase, it was interesting or even necessary to monitor the evolution of $re_m$ during a computation. This, in fact, gives a measure of "reliability" for intermediate results and offers useful indications to adjust the various parameters controlling the process (RTHD, EPS, T, $T_e$, etc.). But, once a particular FP system has been properly setup, associating $re_m$ to every number stored in the system is not necessary for a functioning implementation of the method, and it is practically useless for applications! In fact, if the condition (5.2) is enforced throughout every computation, when one or two operands arrive at the FP processor we know that they must be with $re_m < RTHD$, even if they have no $re_m$ explicitly associated with them. It is thus sufficient to verify that also the current computation if performed with $re < RTHD$, to ensure that the overall computation has been performed, up to the current point, with $re_m < RTHD$.

This fact has some positive practical consequences:

a) the fpe structure reduces to the couple (2.1), instead of the triple (5.1), we can thus save the memory required to store $re_m$;

b) in a specialized FP processor we can save the module dedicated to the computation of $re_m$, i.e. the block marked with "max" in the scheme of Fig. 18.

We must say that, after having focused our proposal quite well, Pirus then missed to classify it in the table resuming the known techniques! Keyser, instead, classified the technique in a category entitled: "Estimated error" [33]. Then Pirus adds some final considerations:

> "The computation time also increases since it is necessary for each operation to evaluate the error. He proposes to reduce this time by using parallel machines. This would make possible to execute several calculations simultaneously, but the accumulation of operations introduced in the calculations should still result in a significant increase in execution time, since the determination of an error is made at each operation of a calculation. For example, the calculation of 3 + 4 * x * x * x - 5 requires 5 elementary operations plus a dozen additional operations, not including the determination of the rounding error.
> In addition, to compare two numbers, one need to calculate and compare two intervals, which also takes time."

Well, this is more or less true, there is no free meal in this world, but, as far as we know, the other alternatives, even those proposed by Pirus further on in her thesis (a combination of rational arithmetic and of the permutation-perturbation method[12]), are even worse in performance terms, plus they suffer of other drawbacks, like lack of generality, insufficient robustness or necessity of extensively rewriting the application code and of specialized programming, as we have seen in section 1.

Moreover, at the hardware level, we can resort both to parallelism and to pipelining of operations, as we have seen in section 11. We believe a lot could be done in this sense and the final throughput of an optimized FP processor such conceived could perhaps approach that of a traditional processor.

Finally, we should mention our recent discovery (in these days!) of [39]. This work shares with ours the same basic idea, but re-proposed as a novelty in recent years! Apparently the authors were unaware of our work [25], as it does not appear in the references; this confirms that our work went mostly neglected, as we wrote above, however, the fact that the same idea has been proposed independently by other researchers, is reassuring.

The paper by Lang and Bruguera is a much lighter work than this one, thus only the basic ideas and some examples are presented. The authors move the first steps along the same path followed here, and then some differences appear. The main one is that, in order to overcome the difficulties of the relative error in vicinity of 0, they used the scaled absolute error, in place of the absolute error and of the non-standard definition of the relative error that we have used here. Pro and cons of this choice should be investigated; currently we do not see reasons to prefer this approach rather than ours. Some examples which show a good behavior of their error estimation with respect to the true error are presented, but no extensive statistical validation of the results is given. In fact, as computations are performed under no constraints, they cannot guarantee for accuracy of the estimation, which, for admission of the authors, is mainly "qualitative". The method, as described, is thus intended primarily for verification of algorithms, not to provide results within given error bounds; no dynamic precision extension is thus advocated. The authors conclude with their expectation to develop a hardware implementation of the method, but no detail is given on this task, which appears quite underestimated.

---

[12] In the method proposed by Pirus each number is represented as a 5-tuple: *(x, sd, num, den, s)*, where *x* is the value represented as a fl.p. number, *sd* is the number of significant digits of *x*, and *(num/den)\*s* is a rational approximation of *x* (*num* and *den* are unsigned integers with arbitrary number of digits and *s* is the sign of the rational); *sd* is evaluated with a permutation-perturbation method, which requires re-executing 3 times the algorithm with perturbed data. In case *sd* for some of the results is too small, then the algorithm is re-executed with exact computations resorting to rational arithmetic.

    Comparing this method to ours it's easy to point out that: a) it's far more complicated; b) the notion of error estimation is substituted with that of "number of significant digits", which is similar, but less precise; c) it requires reprogramming of the applications; d) it requires more resources in terms of space; e) most likely also execution time would be longer, even recurring to concurrent evaluation of the 3 perturbed problems (which implies the use of a truly parallel machine), particularly if frequent resort to rational computations is required; f) the confidence level of *sd* is low, due to the limits of the permutation-perturbation method (see section 1), this lowers the robustness of the method.

## 13. Conclusions

We have proposed a method of providing a realistic error estimation for each floating-point (FP) number resulting from a computation. The method compares favorably with interval arithmetic and other methods of error analysis, which require reprogramming of the applications, give overly pessimistic results or provide insufficient confidence levels for the estimated error, moreover at an unacceptable cost in execution time.

With an appropriate implementation the error estimation has a compact statistical distribution, so that it can be used to determine a narrow confidence interval where the real value of a FP number has an extremely high probability of falling. This probability can be made virtually 1 in a practical system.

The error estimation is also used to monitor the relative error during an entire sequence of computation, thereby ensuring the validity of the results.

Besides getting well conditioned problems solved robustly, the error estimation also enables the *detection* of badly conditioned problems, which would require a higher precision to be solved adequately. In these cases, at a certain point, the maximum relative error will exceed a threshold of acceptability. If such critical conditions should arise, an automatic, *dynamic extension* of precision can take place, and the latest operations can be re-executed with greater precision. It is the availability of the error estimation which provides the opportunity for dynamic extension of precision. In this way, the computations can be constrained not to exceed a maximum relative error; consequently, we have shown that the probability that these calculations are performed within given error bounds can be made extremely high.

Each number is represented by a data structure, called *fpe*, which collects the information on the computed value and on the estimated error. The maximum relative error occurring during a computation is monitored, but it is not necessary to associate it with each number and store it in the structure fpe (which thus reduces to a couple). As the error estimation can be represented with less precision than the value part, a typical format for this structure takes less than twice the space reserved to represent the value of the number, which compares favorably with the space requirements of other methods.

We have outlined the design of a hardware processor implementing the arithmetic operations, and the square root, on objects of type *fpe.* The proposed solution exploits the intrinsic parallelism of the operations involved. Further improvements could be achieved by pipelining the operations. A lot can be done in this sense and the final throughput of an optimized FP processor incorporating the method could perhaps approach that of a traditional processor.

Glauco Masotti was born in Ravenna, Italy, in 1955. He graduated, summa cum laude, Doctor in Electronic Engineering at the University of Bologna, Italy, in 1980. Since then he has been involved in the design and development of various software systems, but he also got patents for some electronic devices. He specialized in CAD/CAM systems, and in research in the field. Most of his work, being proprietary, has not been published. In early years, in Bologna, he cooperated at the development of the GBG drafting system for CAD.LAB (former name of Think3). For COPIMAC he led the development of Aliseo, one of the first highly interactive, parametric and feature-based, solid modelers. In 1989, he was visiting scholar at the University of Southern California (USC), in Los Angeles, where he worked at the design of an object-oriented geometric toolkit. He also proposed an extension to the C++ language, known as EC++ (Extended C++). Back in Italy, he undertook research, at the Dept. of Electronics and Computer Science of the University of Bologna, on symbolic methods and numerical problems in geometric algorithms. He also taught a course on Industrial Automation, with emphasis on geometric modeling and CAD/CAM integration. In 1992 he joined Ecocad (later acquired by Think3), in Pesaro, Italy, where he was in charge of the design and development of various parts of Eureka, an advanced surface and solid modeler. One of his major contributions was the development of a module for assembly modeling and kinematic simulation, supporting the positioning of parts via mating constraints. In 1997 he was co-founder of IeS, Ingegneria e Software Srl, in Bologna, where he worked till 2001, developing software for automatic conversion of drawings in solid models, automatic production of exploded views, interpolation of measured 3D points, plane development of surfaces with minimal distortion, and finite element analysis. He was also involved in structural analysis and CFD problems. In latest years, he undertook autonomous researches in various fields.

He is an enthusiast of sailing, which he practiced at a competitive level until 2002, racing on catamarans.